\documentclass[8pt]{article}

\usepackage{hyperref}
\usepackage{graphicx}
\usepackage{amsmath}
\usepackage{xcolor}
\usepackage{amssymb}
\usepackage{amsfonts}
\usepackage{mathabx}
\usepackage{marginnote}
\usepackage{soul}
\def\ddd{n}

\def\omegen{{\widetilde\Omega}}

\newcommand{\mcD}{\mathcal{D}}

\newcommand{\GGG}{\mathcal{G}}

\newcommand{\mcI}{\mathcal{I}}

\newcommand{\mcL}{\mathcal{L}}

\newcommand{\mbR}{\mathbb{R}}

\newcommand{\mbRn}{{\mathbb{R}^n}}

\def\omgsc{{\omgs\cup\omgc}}
\newcommand{\omg}{{\Omega}}
\newcommand{\omgs}{{\Omega}}
\newcommand{\omgc}{{\Omega_{\mathcal I}}}

\def\rchi{\raise2.2pt\hbox{$\chi$}}
\def\rrchi{\raise1.pt\hbox{$\chi$}}

\newcommand{\beq}{\begin{equation}}
\newcommand{\eeq}{\end{equation}}

\def \alphab{{\boldsymbol\alpha}}

\def \nub{{\boldsymbol \nu}}

\def \bthe{{\boldsymbol \Theta}}

\def \xb{\mathbf{x}}

\def \yb{\mathbf{y}}

\def \wtu {\widetilde u}
\def \wtuN {\widetilde u_N}

\newtheorem{theorem}{Theorem}[section]

\begin{document}

\begin{center}
{\bf \Large The fractional Laplacian operator on \\ bounded domains as a special case of the \\ nonlocal diffusion operator}\\

\vspace*{15pt}
\renewcommand{\thefootnote}{\fnsymbol{footnote}} 	
\small MARTA D'ELIA \footnote{Corresponding author, Phone: +1(404) 345-2587.}
\setcounter{footnote}{0}

{\footnotesize Department of Scientific Computing,
Florida State University, 400 Dirac\\ Science Library, Tallahassee FL 32306-4120, USA; {\tt mdelia@fsu.edu}.}

\vspace*{8pt}

\small MAX GUNZBURGER

{\footnotesize Department of Scientific Computing, Florida State University, 
400 Dirac\\ Science Library, Tallahassee FL 32306-4120, USA; {\tt gunzburg@fsu.edu}.}
\end{center}

\vspace{5pt}

\footnotesize
{\bf Abstract.}
We analyze a nonlocal diffusion operator having as special cases the fractional Laplacian and fractional differential operators that arise in several applications. In our analysis, a nonlocal vector calculus is exploited to define a weak formulation of the nonlocal problem. We demonstrate that, when sufficient conditions on certain kernel functions hold, the solution of the nonlocal equation converges to the solution of the fractional Laplacian equation on bounded domains as the nonlocal interactions become infinite. We also introduce a continuous Galerkin finite element discretization of the nonlocal weak formulation and we derive a priori error estimates. Through several numerical examples we illustrate the theoretical results and we show that by solving the nonlocal problem it is possible to obtain accurate approximations of the solutions of fractional differential equations circumventing the problem of treating infinite-volume constraints.\\



\normalsize
%


\pagestyle{myheadings}
\thispagestyle{plain}
\markboth{THE FRACTIONAL LAPLACIAN AS A NONLOCAL DIFFUSION OPERATOR}{M. D'ELIA, M. GUNZBURGER}

\section{Introduction and motivation}

Nonlocal models have been recently used in many applications, including continuum mechanics \cite{sill:00}, graph theory \cite{lova:06}, nonlocal wave equations \cite{weck:05}, and jump processes \cite{barl:09,bass:10,burc:11}; we consider nonlocal diffusion operators which arise in several and diverse applications such as image analyses \cite{bucm:10,gilboa:595,gilboa:1005,lzob:10}, machine learning \cite{robd:10}, nonlocal Dirichlet forms \cite{appl:04}, kinetic equations \cite{bass:84,limi:10}, phase transitions \cite{bach:99,fife:03}, nonlocal heat conduction~\cite{bodu:09}, and the peridynamic model for mechanics \cite{chen:11,sill:00}. In this work we consider a nonlocal integral operator for {\it anomalous diffusion} which has, as special cases, the fractional Laplacian and fractional derivative operators that are commonly used to model anomalous diffusion \cite{mekl:04}. Physical phenomena exhibiting this property cannot be modeled accurately by the usual advection-dispersion equation; among others, we mention turbulent flows \cite{carr:01, shle:87} and chaotic dynamics of classical conservative systems \cite{zasl:93}.

Nonlocal models differ from the classical partial differential equation models in the fact that in the latter case interactions between two domains occur only due to contact, whereas in the former case interactions can occur at a distance. In particular, let $\Omega\subset\mbRn$ denote a bounded, open domain. For $u(\xb)\colon \Omega \to \mbR$, define the action of the nonlocal diffusion operator $\mcL$ on the function $u(\xb)$ as
\begin{displaymath}
   \mcL u(\xb) := 2\int_{\mbR^n} \big(u(\yb)-u(\xb)\big) \, \gamma (\xb, \yb )\,d\yb \qquad  \forall\,\xb \in \omgs \subseteq \mbRn,
\end{displaymath}
where the volume of $\omgs$ is non-zero and the {\it kernel} $\gamma (\xb, \yb )\colon\Omega\times\Omega\to\mbR$ is a non-negative symmetric mapping. We are interested in the nonlocal, steady-state diffusion equation 
\begin{displaymath}
\left\{\begin{aligned}
-\mcL u = f &\qquad \mbox{on $\omgs$} \\
u  = 0 &\qquad \mbox{on $\omgc$},
\end{aligned}\right.
\end{displaymath} 
where the equality constraint (extension of a Dirichlet boundary condition for differential problems) acts on an interaction volume $\omgc$ that is disjoint from $\omgs$.
 
The numerical solution of fractional differential equations is an open problem and it is the object of current research in applications of models of fractional order; see \cite{podl:09} for recent work including many citations to the literature. Common techniques include methods that take advantage of Laplace and Fourier transforms to obtain classical solutions \cite{lian:06, valk:05} and finite difference methods \cite{memo:04,meta:04,sche:08,shen:04,yust:06} used for constructing numerical approximations. Galerkin discretizations and their error analysis have been considered in \cite{ervi:05,fix:04} for the discretization of the steady state fractional advection-diffusion equation.

A goal of this paper is to develop and analyze discretization methods for fractional Laplacian equations on bounded domains. We do this by exploiting the fact that {\em the fractional Laplacian operator $(-\Delta)^s$ is a special case of the nonlocal operator $\mcL$.} In particular, we compare the solution of the nonlocal steady diffusion equation with the solution of the fractional Laplacian equation on bounded domains and show that solving nonlocal problems is a viable alternative to solving fractional differential equations. A main contribution of this work is to show that not only $(-\Delta)^s$ is a special case of $\mcL$, but that it is the limit of the nonlocal operator as the nonlocal interactions become infinite, provided that sufficient conditions on certain kernel functions hold. This fact has important consequences in the treatment of problems involving the fractional Laplacian equations on bounded domains. In fact, nonlocal problems are a well posed and are a more general formulation of fractional differential models; this is useful for both the analysis of fractional differential equations and for developing finite-dimensional discretization schemes. In \cite{dglz_2012}, finite-dimensional approximations of nonlocal problems is discussed, including finite element discretizations; Galerkin formulations are introduced and results are proved about their well posedness and about estimates for the approximation error and the condition number of finite element matrices. This helps in designing efficient numerical methods for the solution of nonlocal diffusion problems and, as a consequence, of fractional differential problems. Also, being able to quantify the discrepancy between nonlocal solutions and solutions of fractional differential equations (as we demonstrate in this paper) allows us to determine to what extent the former are accurate approximations of the latter. Furthermore, an important advantage of approximating the solution of fractional differential problems with nonlocal solutions is that in the latter case we do not have to deal with infinite-volume constraints as happens for the solution of fractional differential equations on bounded domains where some expedients for constructing finite-dimensional approximations have to be introduced. 

We note that the analysis and approximation of fractional Laplacian problems for the case $\omgs=\mbR^n$, though of interest in many applications, is not a goal of this work and is not discussed. Here, our interest is strictly on bounded domain problems. 

We also note that a significant advantage of recasting fractional Laplacian problems in terms of the nonlocal problems we consider is that, for the first time, such problems can be treated on bounded domains for the case of $s\le1/2$. Indeed, the key to this is our introduction of volume constraints as a generalization of boundary conditions; the latter are not well defined for $s\le1/2$ because traces of functions in the energy space associated with fractional Laplacian operators are themselves not well defined.

In our analysis, a recently developed nonlocal vector calculus \cite{Du10} is exploited to define a weak formulation of the nonlocal problem and to prove the convergence of the nonlocal solution to the solution of the fractional differential problem. In \S\ref{sec:nlvc}, we provide a brief review of those aspects of the nonlocal calculus that are useful in the remainder of the paper, introduce the kernel function, and discuss its properties. In \S\ref{nonlocal_fractional}, we define the fractional Laplacian as a special case of the nonlocal operator $\mcL$ and prove the convergence of the nonlocal operator to the fractional Laplacian as the nonlocal interactions become infinite. In \S\ref{fdapp} we introduce finite-dimensional discretizations of the nonlocal problem and we study the convergence of the approximate nonlocal solutions to the solution of the fractional Laplacian equation. We also discuss the choice of nonuniform grids for mitigating the high computational costs that occur when the extent of nonlocal interactions becomes large. In \S\ref{sec:num_tests}, we present results of some numerical tests for finite element discretizations of one-dimensional problems; by providing qualitative and quantitative comparisons between approximate nonlocal solutions and exact solutions of the fractional Laplacian equation, these results illustrate the theory introduced in \S\S\ref{nonlocal_fractional} and \ref{fdapp}.

\section{Elements of a nonlocal vector calculus}\label{sec:nlvc}

In this section, we review relevant aspects of the nonlocal calculus including nonlocal operators, kernel functions, and nonlocal function spaces. Details about the nonlocal calculus are found in \cite{Du10}. 

The action of the {\em nonlocal divergence operator} $\mcD\colon \mbRn \to \mbR$ on $\nub$ is defined as
\begin{subequations}\label{ndivgrad}
\begin{equation}\label{ndiv}
 \mcD\big(\nub\big)(\xb) := \int_{\mbRn} \big(\nub(\xb,\yb)+\nub(\yb,\xb)\big)\cdot\alphab(\xb,\yb)\,d\yb\qquad
        \mbox{for $\xb\in\mbRn$},
\end{equation}
where $\nub(\xb,\yb), \alphab(\xb,\yb) \colon \mbRn\times\mbRn\to \mbRn$ with $\alphab$ antisymmetric, i.e., $\alphab(\yb,\xb)=-\alphab(\xb,\yb)$, are given mappings. This definition is justified in \cite{Du10}. The action of the {\em nonlocal gradient operator} $\GGG\colon  \mbRn\times\mbRn\to\mbRn$ on $u$ is defined as
\begin{equation}\label{ngra}
\GGG\big(u\big)(\xb,\yb) := \big(u(\yb)-u(\xb)\big)  \alphab(\xb,\yb) \qquad\mbox{for $\xb,\yb\in\mbRn$},
\end{equation}
\end{subequations}
where $u(\xb)\colon \mbRn$ $\to\mbR$ is a given mapping. The fact that $-\GGG$ and $\mcD$ are adjoint operators is shown in \cite{Du10}; in fact, in \cite{Du10}, the operator $\GGG$ is denoted by $-\mcD^\ast$. 

The operator $\mcL\colon \mbRn \to \mbR$ is defined by
\begin{equation}\label{nonlocalL}
\begin{aligned}
-\mcL u (\xb) &:= \mcD\big(\bthe \cdot \GGG u)(\xb)= 2\int_{\mbRn}\big(u(\yb)-u(\xb)\big) \gamma(\xb,\yb) \,d\yb\\ 
&\quad\mbox{for $\,\,\,\xb\in\mbRn\,\,\,$ with $\,\,\,\gamma(\xb,\yb):=\alphab(\xb,\yb)\cdot \big(\bthe(\xb,\yb)\cdot\alphab(\xb,\yb)\big)$,}
\end{aligned}
\end{equation}
i.e., it is the composition of the nonlocal divergence and gradient operators, where $\bthe(\xb,\yb)=\bthe(\yb,\xb)$ denotes a second-order, positive definite symmetric tensor. Note that the {\it kernel} $\gamma(\xb,\yb)$ is symmetric, i.e., $\gamma(\xb,\yb)=\gamma(\yb,\xb)$. The operator $-\mcL$ is non-negative because $\bthe$ is positive definite and $\mcD$ and $-\GGG$ are adjoint operators. If $\bthe$ is the identity tensor, $\mcL$ can be interpreted as a nonlocal Laplacian operator. 

The {\em interaction domain} corresponding to an open subset $\omgs\subset\mbRn$ is defined by 
$$
   \omgc := \{ \yb\in\mbRn\setminus\omgs \quad\mbox{such that}\quad \alphab(\xb,\yb)\ne{\bf 0}\quad \mbox{for some $\xb\in\omgs$}\}.
$$
Thus, $\omgc$ consists of those points outside of $\omgs$ that interact with points in $\omgs$. The situation $\omgc=\mbRn\setminus\omgs$, which we consider here, is allowable as is the case $\omgs=\mbRn$ which we do not consider due to our interest in fractional Laplacian problems on {\em bounded} domains. 

\subsection{The kernel}\label{sec:kernel}

We assume that the domain $\omgs$, as well as $\omgc$ and $\omgsc$, are  bounded with piecewise smooth boundary and satisfy the interior cone condition. For the kernel $\gamma$, we assume that it is symmetric and satisfies
\begin{equation}\label{gamma-conds}
\left\{\begin{aligned}
\gamma(\xb,\yb)  \geq  0 \quad &\forall\, \yb\in B_\lambda(\xb)
\\
\gamma(\xb,\yb) \ge \gamma_0 >  0 \quad &\forall\, \yb\in B_{\lambda/2}(\xb)
\\
\gamma(\xb,\yb)  = 0 \quad &\forall\, \yb\in (\omgsc) \setminus B_\lambda(\xb)
\end{aligned}\right.
\end{equation}
for all $\xb\in\omgsc$, where $\gamma_0$ and $\lambda$ are given positive constants and $B_\lambda({\xb}) := \{ \yb \in\omgsc  \colon |\yb-\xb|\le \lambda \}$; thus, nonlocal interactions are limited to a ball of radius $\lambda$ which itself is referred to as the {\em interaction radius}. Note that then
\begin{equation}\label{omgie}
  \omgc = \{ \yb\in \mbR^\ddd\setminus\omg \,\,\,\colon\,\,\, |\yb-\xb|<\lambda \mbox{ for $\xb\in\omg$}\},
\end{equation}
i.e., the interaction domain $\omgc$ is a layer of thickness $\lambda$ surrounding $\omg$. We further assume that
\begin{equation}\label{case1}
\frac{\gamma_1}{|\yb-\xb|^{\ddd+2s}} \leq \gamma(\xb,\yb) \leq \frac{\gamma_2}{|\yb-\xb|^{\ddd+2s}} \qquad \mbox{for $\yb\in B_\lambda({\xb})$} 
\end{equation}
for all $\xb\in\omgs$, where $s\in (0,1)$ and $\gamma_1$ and $\gamma_2$ denote positive constants. An example is given by
\begin{equation}\label{kernel_case1}
   \gamma(\xb,\yb) = \frac{\sigma(\xb,\yb)}{|\yb-\xb|^{\ddd+2s}}
\end{equation}
with $\sigma(\xb,\yb)$ bounded from above and below by positive constants.

In \cite{dglz_2012}, other choices for the kernel are discussed and analyzed, showing that the nonlocal diffusion operator $\mcL$ has more general applicability than to just fractional Laplacian problems.

Note that we allow for all $s\in(0,1)$ and not just $s>1/2$.

\subsection{Equivalence of spaces}\label{sec:equivsp}

We respectively define the nonlocal energy semi-norm, nonlocal energy space, and nonlocal volume-constrained energy space by
\begin{subequations}
\begin{align}
&|||v|||   := 
\bigg(\frac12\int_\omgsc\int_{\omgsc}\GGG(v)(\xb,\yb )\,\cdot\big(\bthe(\xb,\yb )\cdot\GGG(v)(\xb,\yb )\big)\,d\yb \, d\xb\bigg)^{1/2} \label{energynorm}
\\
&V(\omgsc)  := \left\{ v  \in L^2(\omgsc) \,\,:\,\, |||v||| < \infty \right\} \label{vspace}\\
&V_c(\omgsc)  := \left\{v\in V(\omgsc) \,\,:\,\, v=0\;{\rm on}\;\omgc\right\}.\label{vcspace}
\end{align}
\end{subequations}
In \cite{dglz_2012}, it is shown that, for kernels satisfying \eqref{gamma-conds} and \eqref{case1}, the nonlocal energy space $V(\omgsc)$ is equivalent to the fractional-order Sobolev space $H^s(\omgsc)$.\footnote{For $s\in(0,1)$ and for a general domain $\omegen\in\mbR^\ddd$, let
$$
|v|_{H^s(\omegen)}^2 :=
\int_\omegen\int_\omegen\frac{\big(v(\yb)-v(\xb)\big)^2}{|\yb-\xb|^{\ddd+2s}}\,d\yb d\xb.
$$
Then, the space $H^s(\omegen)$ is defined by \cite{Adams}
$
H^s(\omegen) := \left\{v\in L^2(\omegen) : \|v\|_{L^2(\omegen)} +
|v|_{H^s(\omegen)}<\infty\right\}.
$
}
This implies that $V_c(\omgsc)$ is a Hilbert space equipped with the norm $|||\cdot|||$. In particular, we have
\begin{equation}\label{eq:equivalence}
C_1\|v\|_{H^s(\omgsc)}\leq |||v||| \leq C_2 \|v\|_{H^s(\omgsc)} \;\;\forall\,v\in V_c(\omgsc)
\end{equation}
for some positive constants $C_1$ and $C_2$. As a consequence, any result obtained below involving the energy norm $|||\cdot|||$ can be immediately reinterpreted as a result involving the norm $\|\cdot\|_{H^s(\omgsc)}$.

\section{Relations between the nonlocal Laplacian and the fractional Laplacian}\label{nonlocal_fractional}

In this section, we introduce the fractional Laplacian operator as a special case of the nonlocal operator $\mcL$. Then, we show that for the {\em special} kernels given by \eqref{kernel_case1} with $\sigma=$ constant, the nonlocal operator $\mcL$ approaches the fractional Laplacian as the extent of nonlocal interactions increases, i.e., as $\lambda\to\infty$.

The fractional Laplacian is the pseudo-differential operator with Fourier symbol $\mathcal{F}$ satisfying \cite{appl:04}
\begin{equation}\label{eq:fourier_fl}
\mathcal{F}\big((-\Delta)^{s}u\big)(\xi) = |\xi|^{2s}\widehat{u}(\xi), \qquad 0< s \leq 1,
\end{equation}
where $\widehat{u}$ denotes the Fourier transform of $u$. Using Fourier transforms, it can be shown \cite{appl:04} that an equivalent characterization of the fractional Laplacian is given by 
\begin{equation}\label{eq:fl}
(-\Delta)^{s}u = c_{\ddd,s} \int_{\mbRn} \frac{u(\xb)-u(\yb)}{|\yb - \xb|^{\ddd+2s}} \, d\yb, \qquad 0< s <1,
\end{equation}
where the normalizing constant $c_{\ddd,s}$ is defined as
\begin{equation}\label{eq:csn}
c_{n,s} = s\,2^{2s}\dfrac{\Gamma\left(\frac{n+2}{2}\right)}{\Gamma\left(\frac{1}{2}\right)\Gamma(1-s)}.
\end{equation}
Here, $\Gamma(\cdot)$ is defined via the improper integral
$$
\Gamma(a)=\int_0^\infty t^{a-1}e^{-t}\,dt
$$
for all complex numbers $a$ having positive real part. Thus, if the kernel $\gamma (\xb, \yb )$ in \eqref{nonlocalL} is defined as 
\begin{equation}\label{kernell}
  \gamma (\xb, \yb ) =  \frac{c_{\ddd,s}}{2|\yb - \xb|^{n+2s}}\qquad\forall\, \xb,\yb\in\mbR^n,
\end{equation}
then
\begin{equation}\label{factor}
   -\mcL = \big(-\Delta\big)^{s}, \qquad 0< s <1,
\end{equation}
establishing that the fractional Laplacian is a special case of the operator $\mcL$ for the choice of $\gamma(\xb,\yb)$ proportional to $1/|\yb - \xb|^{\ddd+2s}$. 

We consider the volume constrained problem
\begin{equation}\label{eq:fl_strong}
\left\{\begin{aligned}
-\mcL u = f & \quad {\rm in} \; \omgs \\
u = 0 & \quad {\rm in} \; \mbR^n\backslash\omgs,
\end{aligned}\right.
\end{equation}
where $\mcL$ is defined as in \eqref{nonlocalL} with $\gamma(\xb,\yb)$ given by \eqref{kernell} and $f\in L^2(\omgs)$. A weak formulation of \eqref{eq:fl_strong} is 
\begin{equation}\label{eq:fl_weak}
\frac{c_{\ddd,s}}{2}\int_{\mbR^n}\int_{\mbR^n}  \frac{u(\yb)-u(\xb)}{|\xb-\yb|^{\ddd+2s}} \big(v(\yb)-v(\xb)\big)\, d\yb\,d\xb = \int_{\omgs} f\,v \, d \xb \quad \forall \; v\in H^s_\omgs(\mbRn),
\end{equation}
where $H^s_\omgs(\mbRn):=\{w\in H^s(\mbR^n): \, w = 0 \; {\rm in} \;\, \mbR^n\backslash \omgs \} $. The existence and uniqueness of such $u\in H^s_\omgs(\mbRn)$ for all $s\in(0,1)$ is proven in \cite{webb}. 

Next, for all $\xb\in\mbRn$, we define
\begin{equation}\label{eq:eps_kernel}
\widetilde\gamma(\xb,\yb) = \left\{
\begin{array}{ll}
\dfrac{ c_{\ddd,s}}{2|\xb-\yb|^{\ddd+2s}} & \yb\in B_\lambda(\xb)\\
[3mm]
0& \yb\in \mbR^n\backslash B_\lambda(\xb).
\end{array}\right.
\end{equation}
Letting $\widetilde\mcL$ denote the nonlocal operator associated with $\widetilde\gamma$, we consider the problem
\begin{equation}\label{eq:eps_strong}
\left\{\begin{aligned}
-\widetilde\mcL \,\wtu = f & \quad {\rm in} \; \omgs \\
\wtu = 0 & \quad {\rm in} \; \omgc,
\end{aligned}\right.
\end{equation}
where, again, the interaction domain $\omgc$ is defined as $\omgc = \{\yb\in\mbR^n\backslash\omgs: \; |\xb-\yb|\leq \lambda \,\, \forall \, \xb\in\omgs \}$. A weak formulation of \eqref{eq:eps_strong} is 
\begin{equation}\label{eq:eps_weak}
\begin{aligned}
\frac{c_{\ddd,s}}{2} \int_{\omgsc}\int_{(\omgsc)\cap B_\lambda(\xb)} \hspace{-2pt} \frac{\wtu(\yb)-\wtu(\xb)}{|\xb-\yb|^{\ddd+2s}}& \big(v(\yb)-v(\xb)\big)\, d\yb\,d\xb \\&= \int_{\omgs} \hspace{-2pt} f\,v \, d \xb \;\;\; \forall \; v\in V_c(\omgsc)
\end{aligned}
\end{equation}
where $V_c(\omgsc)$, isomorphic to $H^s_\omgs(\mbR^n)$ \cite{webb}, is the volume constrained energy space associated with the kernel $\widetilde \gamma$. 

We now show that the solution $u$ of \eqref{eq:fl_weak} is the limit, as $\lambda\to\infty$, of the solution $\wtu$ of \eqref{eq:eps_weak}.

\begin{theorem}
Let $u\in H^s_\omgs(\mbRn)$ and $\wtu\in V_c(\omgsc)$ denote the solutions of \eqref{eq:fl_weak} and \eqref{eq:eps_weak}, respectively. Then, with $I:=\min\{R: \omg\subset B_R(\xb)\; \forall\,\xb\in\omg\}$,
\begin{equation}\label{eq:lambda_limit}
\begin{array}{l}
\|u-\wtu\|_{H^s(\omgsc)}\leq \dfrac{K_n}{C_1^2 s (\lambda-I)^{2s}}\; \|u\|_{L^2(\omgs)}\\
[5mm]
\|u-\wtu\|_{L^2(\omg)}\leq \dfrac{K_n}{C_1^2 s (\lambda-I)^{2s}}\; \|u\|_{L^2(\omgs)},
\end{array}
\end{equation}
where $C_1$ is the equivalence constant in \eqref{eq:equivalence} and $K_n$, independent of $s$ and $\lambda$, is determined by the spatial dimension $n$.
\end{theorem}

{\it Proof.}
Combining \eqref{eq:fl_weak} and \eqref{eq:eps_weak} we have 
$$
\begin{aligned}
 \int_{\mbR^n}\int_{\mbR^n}& \frac{u(\yb)-u(\xb)}{|\xb-\yb|^{\ddd+2s}} \left(v(\yb)-v(\xb)\right)\, d\yb\,d\xb\,-\\ 
 \int_{\omgsc}\int_{(\omgsc)\cap B_\lambda(\xb)}& \frac{\wtu(\yb)-\wtu(\xb)}{|\xb-\yb|^{\ddd+2s}} \left(v(\yb)-v(\xb)\right)\, d\yb\,d\xb = 0 \quad
\forall \, v\in V_c(\omgsc).
\end{aligned}
$$
If we define $\omgs_c=\mbR^n\backslash( \omgsc)$, this equality is equivalent to
$$
\begin{aligned}
 \int_{\omgsc}\int_{\omgs_c} \frac{u(\yb)-u(\xb)}{|\xb-\yb|^{\ddd+2s}} \left(v(\yb)-v(\xb)\right)\, d\yb\,d\xb+\\ 
  \int_{\omgsc}\int_{(\omgsc)\cap B_\lambda(\xb)} \frac{u(\yb)-u(\xb)}{|\xb-\yb|^{\ddd+2s}} \left(v(\yb)-v(\xb)\right)\, d\yb\,d\xb+\\
  \int_{\omgsc}\int_{(\omgsc)\backslash B_\lambda(\xb)} \frac{u(\yb)-u(\xb)}{|\xb-\yb|^{\ddd+2s}} \left(v(\yb)-v(\xb)\right)\, d\yb\,d\xb+\\ 
   \int_{\omgs_c}\int_{\mbR^n} \frac{u(\yb)-u(\xb)}{|\xb-\yb|^{\ddd+2s}} \left(v(\yb)-v(\xb)\right)\, d\yb\,d\xb+\\ 
-\int_{\omgsc}\int_{(\omgsc)\cap B_\lambda(\xb)} \frac{\wtu(\yb)-\wtu(\xb)}{|\xb-\yb|^{\ddd+2s}} \left(v(\yb)-v(\xb)\right)\, d\yb\,d\xb =
\end{aligned}
$$
$$
\begin{aligned}
\int_{\omgsc}\int_{(\omgsc)\cap B_\lambda(\xb)} \frac{(u(\yb)-\wtu(\yb))-(u(\xb)-\wtu(\xb))}{|\xb-\yb|^{\ddd+2s}} \left(v(\yb)-v(\xb)\right)\, d\yb\,d\xb+\\ 
\int_{\omgsc}\int_{(\omgsc)\backslash B_\lambda(\xb)} \frac{u(\yb)-u(\xb)}{|\xb-\yb|^{\ddd+2s}} \left(v(\yb)-v(\xb)\right)\, d\yb\,d\xb+\\
\int_{\omgsc}\int_{\omgs_c} \frac{u(\yb)-u(\xb)}{|\xb-\yb|^{\ddd+2s}} \left(v(\yb)-v(\xb)\right)\, d\yb\,d\xb+\\
\int_{\omgs_c}\int_{\mbR^n} \frac{u(\yb)-u(\xb)}{|\xb-\yb|^{\ddd+2s}} \left(v(\yb)-v(\xb)\right)\, d\yb\,d\xb=0.\\ 
\end{aligned}
$$
Thus, 
\begin{equation}\label{eq:energy1}
\begin{aligned}
\left|\int_{\omgsc}\int_{(\omgsc)\cap B_\lambda(\xb)} \frac{(u(\yb)-\wtu(\yb))-(u(\xb)-\wtu(\xb))}{|\xb-\yb|^{\ddd+2s}} \left(v(\yb)-v(\xb)\right)\, d\yb\,d\xb \right|\leq &\\ 
\left|\int_{\omgsc}\int_{(\omgsc)\backslash B_\lambda(\xb)} \frac{u(\yb)-u(\xb)}{|\xb-\yb|^{\ddd+2s}} \left(v(\yb)-v(\xb)\right)\, d\yb\,d\xb\right|+ & \\
\left|\int_{\omgsc}\int_{\omgs_c} \frac{u(\yb)-u(\xb)}{|\xb-\yb|^{\ddd+2s}} \left(v(\yb)-v(\xb)\right)\, d\yb\,d\xb\right|  + &\\
\left|\int_{\omgs_c}\int_{\mbR^n} \frac{u(\yb)-u(\xb)}{|\xb-\yb|^{\ddd+2s}} \left(v(\yb)-v(\xb)\right)\, d\yb\,d\xb\right|= & \\
\mcI_1 + \mcI_2 + \mcI_3.\;\;\; &
\end{aligned}
\end{equation}
We treat $\mcI_1$, $\mcI_2$, and $\mcI_3$ separately. We first split $\mcI_1$ in the two integrals $\mcI_a$ and $\mcI_b$:
$$
\begin{aligned}
\mcI_1 =	 & \left|\int_{\omgsc}\int_{(\omgsc)\backslash B_\lambda(\xb)} \frac{u(\yb)-u(\xb)}{|\xb-\yb|^{\ddd+2s}} \left(v(\yb)-v(\xb)\right)\, d\yb\,d\xb\right|= & \\
& \quad\;\;\left|\int_{\omg}\int_{(\omgsc)\backslash B_\lambda(\xb)} \frac{u(\yb)-u(\xb)}{|\xb-\yb|^{\ddd+2s}} \left(v(\yb)-v(\xb)\right)\, d\yb\,d\xb\right|  & \\
& \;\;\,\;\,\left|\int_{\omgc}\int_{(\omgsc)\backslash B_\lambda(\xb)} \frac{u(\yb)-u(\xb)}{|\xb-\yb|^{\ddd+2s}} \left(v(\yb)-v(\xb)\right)\, d\yb\,d\xb\right|= & \mcI_a + \mcI_b.
\end{aligned}
$$
Then,
$$
\begin{aligned}
\mcI_a = & \left|\int_{\omg} u(\xb)v(\xb)\int_{(\omgsc)\backslash B_\lambda(\xb)} \frac{1}{|\xb-\yb|^{\ddd+2s}} \, d\yb\,d\xb\right| \leq  \\
& \left|\int_{\omg} u(\xb)v(\xb)\int_{B_{\lambda+I}(\xb)\backslash B_\lambda(\xb)} \frac{1}{|\xb-\yb|^{\ddd+2s}} \, d\yb\,d\xb\right|.   \\
\end{aligned}
$$
The value of the inner integral, say $i$, depends on the dimension $n$ of the problem:
\begin{equation}\label{constants}
n=1, \;\; i=\frac{2}{s\lambda^{2s}}; \qquad n=2, \;\; i=\frac{2\pi}{s\lambda^{2s}}; \qquad n=3, \;\; i=\frac{4\pi}{s\lambda^{2s}}. 
\end{equation}
With $k_n= 2,2\pi,4\pi$ for $n=1,2,3$ respectively, we thus have
\begin{equation}\label{Ia}
\begin{aligned}
\mcI_a & \leq \frac{k_n}{s\lambda^{2s}} \left|\int_{\omg} u(\xb)v(\xb)\,d\xb\right|\leq 
\frac{k_n}{s\lambda^{2s}}\|u\|_{L^2(\omgs)}\|v\|_{L^2(\omgs)}\\
&\leq \frac{k_n}{s\lambda^{2s}}\|u\|_{L^2(\omgs)}\|v\|_{H^s{\omgsc}}\leq
\frac{k_n}{C_1 s\lambda^{2s}}\|u\|_{L^2(\omgs)}|||v|||,
\end{aligned}
\end{equation}
where we used the Cauchy--Schwarz inequality and \eqref{eq:equivalence}. Next,
\begin{displaymath}
\begin{aligned}
\mcI_b = & \left|\int_{\omgc}\int_{(\omgsc)\backslash B_\lambda(\xb)} \frac{ u(\yb)v(\yb)}{|\xb-\yb|^{\ddd+2s}} \, d\yb\,d\xb\right| \leq  \\
& \left|\int_{B_{\lambda+I}(\widehat{\xb})\backslash B_{\lambda-I}(\widehat{\xb})}\int_\omg \frac{ u(\yb)v(\yb)}{|\xb-\yb|^{\ddd+2s}} \, d\yb\,d\xb\right| = \\
& \left|\int_\omg u(\yb)v(\yb) \int_{B_{\lambda+I}(\widehat{\xb})\backslash B_{\lambda-I}(\widehat{\xb})}  \frac{1}{|\xb-\yb|^{\ddd+2s}} \, d\xb\,d\yb\right|,
\end{aligned}
\end{displaymath}
where $\widehat{\xb}$ can be any point in $\omgs$; this implies that we can choose $\widehat{\xb} = \yb$ for any $\yb\in\omgs$. Thus, we have
\begin{equation}\label{Ib}
\begin{aligned}
\mcI_b = & \left|\int_\omg u(\yb)v(\yb) \int_{B_{\lambda+I}(\yb)\backslash B_{\lambda-I}(\yb)}  \frac{1}{|\xb-\yb|^{\ddd+2s}} \, d\xb\,d\yb\right|\leq\\
& \frac{k_n}{s(\lambda-I)^{2s}} \left|\int_{\omg} u(\yb)v(\yb)\,d\yb\right|\leq 
\frac{k_n}{C_1 s(\lambda-I)^{2s}}\|u\|_{L^2(\omgs)}|||v|||,
\end{aligned}
\end{equation}
where the value of the inner integral is obtained in a way similar to \eqref{constants}. Next, we treat $\mcI_2$ and $\mcI_3$.
\begin{equation}\label{I2}
\begin{aligned}
\mcI_2 \leq & \left|\int_{\omg} u(\xb)v(\xb)\int_{\mbR^\ddd\backslash B_\lambda(\xb)} \frac{1}{|\xb-\yb|^{\ddd+2s}} \, d\yb\,d\xb\right|\leq \\
& \frac{k_n}{s\lambda^{2s}} \left|\int_{\omg} u(\xb)v(\xb)\,d\xb\right|\leq 
\frac{k_n}{C_1 s\lambda^{2s}}\|u\|_{L^2(\omgs)}|||v|||.
\end{aligned}
\end{equation}
Again, the value of the inner integral is obtained in a way similar to \eqref{constants}. The same bound can be found for $\mcI_3$. Now, letting $v=u-\widetilde u$ and combining \eqref{eq:energy1}, \eqref{Ia}, \eqref{Ib}, and \eqref{I2}, we have
$$
|||u-\widetilde u |||^2\leq \frac{4 k_n}{C_1 s(\lambda-I)^{2s}}\|u\|_{L^2(\omgs)} |||u-\widetilde u |||.
$$
Letting $K_n=4k_n$ we conclude that
$$
|||u-\widetilde u |||\leq \frac{K_n}{C_1 s(\lambda-I)^{2s}}\|u\|_{L^2(\omgs)}.
$$
The bound for the $H^s$ norm is obtained using the equivalence relation \eqref{eq:equivalence}; the same bound holds also for the $L^2$ norm because $\|w\|_{L^2(\omgs)}\leq\|w\|_{H^s(\omgsc)}$ for all $w\in H^s(\omgsc)$. $\boxvoid$

According to this theorem, lower values of $s$ lead to slower convergence, i.e., the effects of nonlocality become more pronounced in the sense that the value of the interaction radius $\lambda$ has to be bigger in order to achieve the same error for smaller values of $s$ compared to larger ones.

It is tempting to think that the $L^2$ error estimate in \eqref{eq:lambda_limit} is pessimistic and that a refined analysis would show that that error has a better rate of convergence compared to the $H^s$ error. However, the numerical illustrations of \S\ref{eq:lambda_limit} show that the $L^2$ error estimate in \eqref{eq:lambda_limit} is indeed sharp.   

\section{Finite-dimensional approximations}\label{fdapp}

In this section, we consider the convergence of solutions of finite-dimensional discretizations (including finite element approximations) of the nonlocal problem to solutions of the fractional Laplacian equation. We also discuss the use of nonuniform grids for mitigating the computational costs caused by the increase in the size of the interaction domain as $\lambda\to\infty$. 

We consider conforming {finite-dimensional} subspaces
\begin{equation}\label{conf}
     V^N(\omgsc) \subset V(\omgsc)
\end{equation}
parametrized by an integer $N\to\infty$ and then define the constrained finite-dimensional subspace
\begin{equation}\label{confc}
\begin{aligned}
     V^N_c(\omgsc) &:= V^N(\omgsc)\cap V_c(\omgsc)
     \\&= \left\{v\in V^N(\omgsc) \,\,:\,\, v=0\;{\rm on}\;\omgc\right\}.
\end{aligned}
\end{equation}
The common choice for $N$ is the dimension of the subspaces $V^N(\omgsc)$. We assume that, for any function $v\in V(\omgsc)$, the sequence of best approximations with respect to the energy norm $|||\cdot|||$ converges, i.e., 
\begin{equation}\label{baerror}
  \lim\limits_{N\to\infty} \,  \inf\limits_{v_N\in V^N(\omgsc)} ||| v - v_N||| = 0
  \qquad\forall\, v \in V(\omgsc).
\end{equation}
We also define the Galerkin approximation $\widetilde u_N\in V^N_c(\omgsc)$ of the solution $\widetilde{u}$ of \eqref{eq:eps_weak} as the solution of
\begin{equation}\label{wosd}
\begin{aligned}
\frac{c_{\ddd,s}}{2}\int_{\omgsc}\int_{(\omgsc)\cap B_{\lambda}(\xb)} \frac{\widetilde u_N(\yb)-\widetilde u_N(\xb)}{|\yb-\xb|^{\ddd+2s}} \big(v(\yb)-v(\xb)\big)\,d\yb d\xb = \int_{\omgs} f v\,d\xb,\\
\forall \, v\in V_c^N(\omgsc).
\end{aligned}
\end{equation}
In \cite{dglz_2012} it is shown that $\widetilde  u_N\to \widetilde u$ for $N\to\infty$ for all $f\in L^2(\omgs)$ and $s\in(0,1)$. 

\subsection{Finite element approximations} 

We consider finite element approximations for the case that both $\omgsc$ and $\omgs$ are polyhedral domains. We partition $\omgsc$ into finite elements and denote by $h$ the diameter of the largest element in the partition. We assume that the interface $\overline\omgs\cap{\overline\omg}_{\mathcal I}$ consists of finite element faces and that the partition is shape-regular and quasi-uniform \cite{brsc:08} as the grid size $h\to 0$. We choose the subspace $V_c^N(\omgsc)\subset V_c(\omgsc)$ to consist of piecewise polynomials of degree no more than $m$ defined with respect to the partition. We have that $N\to\infty$ as $h\to0$, where $N$ denotes the dimension of the approximating space $V_c^N(\omgsc)$. 

In \cite{dglz_2012} it is shown that, assuming that $\widetilde u\in H^{m+t}(\omgsc)$ and $s\leq t\leq 1$, for the kernel given in \eqref{kernel_case1} and for sufficiently small $h$, there exists a constant $C_3$ such that
\begin{equation}\label{fe_convergence}
\|\widetilde u - \widetilde u_N\|_{H^s(\omgsc)} \leq C_3 h^{m+t-s} \|\widetilde u\|_{H^{m+t}(\omgsc)}.
\end{equation}
Our goal is to find an estimate of the error $\|u-\widetilde u_{N,\lambda}\|_{H^s(\omgsc)}$, i.e., of the difference between the solution of the fractional Laplacian equation and the finite element solution of \eqref{wosd} for a particular value of the interaction radius $\lambda$.  The estimate is easily found combining \eqref{eq:lambda_limit} and \eqref{fe_convergence}. By the triangle inequality, we have 
\begin{equation}\label{u_uNlambda_estimate}
\begin{aligned}
\|u-\widetilde u_{N,\lambda}\|_{H^s(\omgsc)} \leq &\; \|u-\widetilde u\|_{H^s(\omgsc)} + \|\widetilde u-\widetilde u_{N,\lambda}\|_{H^s(\omgsc)}  \\
\leq & \;\frac{C_4}{s(\lambda-I)^{2s}} \|u\|_{L^2(\omgs)} + C_3 h^{m+t-s} \|\widetilde u\|_{H^{m+t}(\omgsc)},
\end{aligned}
\end{equation}
where $C_4=K_n/C_1^2$. Note that the convergence to the solution of the fractional Laplacian equation depends on both $\lambda$ and $N$; for very fine grids, i.e., very large $N$, the finite element approximation error is negligible whereas, for very large interaction domains, i.e., very large $\lambda$, the finite element approximation error dominates. 

Note also that this result holds for a uniform grid on $\omgsc$, but, in practice, for very fine uniform grids and very large interaction radii, computing $\widetilde u_{N,\lambda}$ is not affordable. However, the contributions of the integrals in \eqref{wosd} over regions of $\omgc$ far from the domain $\Omega$ are not as significant as those of regions close to it; for this reason, it might be convenient to utilize a coarser grid as we move further away from $\omg$.

\subsection{Nonuniform grid in $\omgc$}\label{sec:non_unif}

Let the partition of $\omgsc$ be uniform in $\omgs$ and nonuniform in $\omgc$ and let problem \eqref{wosd} be further discretized by approximating the integrals via numerical integration;\footnote{Problem \eqref{wosd} has also to be solved using numerical integrations for the case of uniform grids, in contrast to case of the Laplacian operator for which the stiffness matrix can be evaluated exactly. However, for the nonlocal case and for uniform grids, the integration error can easily be made negligible whereas nonuniform coarser grids in $\omgc$ may affect the accuracy of the approximation.} we denote by $\wtuN$ the corresponding solution and find a bound for the error $\|\widetilde u - \wtuN\|_{H^s(\omgsc)}$.

Let $A(\cdot,\cdot):V\times V\rightarrow \mbR$ and $F(\cdot):V\rightarrow\mbR$ respectively denote the bilinear form and functional associated with problem \eqref{eq:eps_weak}:
\begin{subequations}\label{forms}
\begin{align}
\label{bilinear} A(w,v) = & \frac{c_{\ddd,s}}{2}\int_{\omgsc}\int_{(\omgsc)\cap B_{\lambda}(\xb)} \frac{w(\yb)-w(\xb)}{|\yb-\xb|^{\ddd+2s}} \big(v(\yb)-v(\xb)\big)\,d\yb d\xb   \\
F(v) = & \int_{\omgs} f v\,d\xb \qquad {\rm for}\; w,v\in V.
\end{align}
\end{subequations}
Let $A_N(\cdot,\cdot):V^N\times V^N\rightarrow \mbR$ denote an approximation of $A$ obtained by numerical integration of \eqref{bilinear} for $w,v\in V^N$ and let $\wtuN$ denote the solution of\footnote{Of course, in practice, one also has to use numerical integration to approximate the right-hand side functional $F(v)$. However, this term does not pose any problems so that we assume it can be integrated as accurately as one needs for it not to affect the accuracy of approximate solutions.}  
\begin{equation}\label{wosdN}
A_N(\wtuN,v_N) = F(v_N) \quad \forall \, v_N\in V^N. 
\end{equation}
We assume that $A_N$ is a coercive (with coercivity constant $C_5$) and continuous bilinear form on $V^N$, then, problem \eqref{wosdN} is well-posed; furthermore, we can apply the Strang lemma for generalized Galerkin problems to obtain an a priori error estimate for the energy norm \cite{ern}. We have 
\begin{equation}\label{strang}
\begin{aligned}
|||\wtu-\wtuN||| & \leq\inf\limits_{v_N\in V^N} \left\{ \left(1+\frac{1}{C_5}\right)|||\wtu - v_N||| \right.\\ 
&\left. +\frac{1}{C_5} \sup_{w_N\in V^N\backslash\{\bf 0\}} \frac{|A(v_N,w_N)-A_N(v_N,w_N)|}{||| w_N|||}  \right\}.
\end{aligned}
\end{equation} 
Next, we find a similar result for the $H^s$ norm. Let $v_N$ be any function in $V^N(\omgsc)$, we first find a bound for the difference $\|\wtuN-v_N\|_{H^s(\omgsc)}$. We have
\begin{displaymath}
\begin{aligned}
C_5 |||\wtuN-v_N|||^2 & \leq A_N(\wtuN-v_N,\wtuN-v_N)\\
& = A_N(\wtuN,\wtuN-v_N) - A_N(v_N,\wtuN-v_N)\\
& = F(\wtuN-v_N) - A_N(v_N,\wtuN-v_N)\\
& = A(\wtu,\wtuN-v_N) - A_N(v_N,\wtuN-v_N)\\
& = A(\wtu-v_N,\wtuN-v_N) + A(v_N,\wtuN-v_N) - A_N(v_N,\wtuN-v_N)\\
& \leq |||\wtu-v_N|||\,|||\wtuN-v_N||| + \left|A(v_N,\wtuN-v_N) - A_N(v_N,\wtuN-v_N)  \right|.
\end{aligned}
\end{displaymath}
Letting $\Delta A(v_N,w_N)= |A(v_N,\wtuN-v_N) - A_N(v_N,\wtuN-v_N)|$, the previous inequality implies that
\begin{displaymath}
\begin{aligned}
|||\wtuN-v_N||| \leq \frac{1}{C_5}\left(|||\wtu-v_N|||+ \frac{\Delta A(v_N,\wtuN-v_N)}{|||\wtuN-v_N|||}\right).
\end{aligned}
\end{displaymath}
Using the equivalence of norms \eqref{eq:equivalence} we have
\begin{displaymath}
\begin{aligned}
&C_1 \|\wtuN-v_N\|_{H^s(\omgsc)} \leq |||\wtuN-v_N ||| \leq \\
&\frac{C_2}{C_5} \|\wtu-v_N\|_{H^s(\omgsc)} + \frac{1}{C_5}\sup_{w_N\in V^N\backslash\{\bf 0\}} \frac{\Delta A(v_N,w_N)}{|||w_N|||}.
\end{aligned}
\end{displaymath}
Thus,
$$
\|\wtuN-v_N\|_{H^s(\omgsc)} \leq \frac{C_2}{C_1 C_5} \|\wtu-v_N\|_{H^s(\omgsc)} + \frac{1}{C_1 C_5}\sup_{w_N\in V^N\backslash\{\bf 0\}} \frac{\Delta A(v_N,w_N)}{|||w_N|||}.
$$
Now, by the triangle inequality, for all $v_N\in V^N(\omgsc)$, we have
\begin{displaymath}
\begin{aligned}
\|\wtu - \wtuN\|&_{H^s(\omgsc)} \leq \|\wtu - v_N\|_{H^s(\omgsc)} + \|\wtuN-v_N\|_{H^s(\omgsc)} \\
& \leq \left(1+\frac{C_2}{C_1 C_5}\right)\|\wtu - v_N\|_{H^s(\omgsc)} + \frac{1}{C_1 C_5}\sup_{w_N\in V^N\backslash\{\bf 0\}} \frac{\Delta A(v_N,w_N)}{|||w_N|||}. \\
\end{aligned}
\end{displaymath}
Next, we consider a finite element discretization; we denote by $\{E_i\}_{i=1}^{N+K}$, $\{E_i\}_{i=1}^{N}\subset\omgs$ and $\{E_i\}_{i=N+1}^{N+K}\subset\omgc$, the elements of the partition, and by $h_i$ the diameter of each element; we also denote the (uniform) grid size in $\omgs$ by $\widehat h = h_i$, for $i=1,\ldots N$.
Note that
\begin{equation}
\|\wtu - v_N\|_{H^s(\omgsc)} = \left\|\sum_{i=1}^{N+K} \left.(\wtu - v_N)\right|_{E_i}\right\|_{H^s(\omgsc)} \leq \sum_{i=1}^{N+K} \| (\wtu - v_N)\|_{H^s(E_i)}.
\end{equation}
We conclude that 
\begin{displaymath}
\begin{aligned}
\|\wtu - \wtuN\|_{H^s(\omgsc)}& \leq \left(1+\frac{C_2}{C_1 C_5}\right) \sum_{i=1}^{N+K} \| (\wtu - v_N)\|_{H^s(E_i)} \\
&+ \frac{1}{C_1 C_5}\sup_{w_N\in V^N\backslash\{\bf 0\}} \frac{\Delta A(v_N,w_N)}{|||w_N|||}. 
\end{aligned}
\end{displaymath}
Because the previous inequality holds for all $v_n\in V^N$, we can write 
\begin{equation}\label{HsStrang}
\begin{aligned}
\|&\wtu - \wtuN\|  _{H^s(\omgsc)}\leq\\
&\leq \inf_{v_N\in V^N} \left\{\left(1+\frac{C_2}{C_1 C_5}\right) \sum_{i=1}^{N+K} \| (\wtu - v_N)\|_{H^s(E_i)} + \frac{1}{C_1 C_5}\sup_{w_N\in V^N\backslash\{\bf 0\}} \frac{\Delta A(v_N,w_N)}{|||w_N|||} \right\} \\
&\leq \left(1+\frac{C_2}{C_1 C_5}\right) \sum_{i=1}^{N} h_i^{m+t-s}\|\wtu\|_{H^s(E_i)} + \frac{1}{C_1 C_5}\inf_{v_N\in V^N} \sup_{w_N\in V^N\backslash\{\bf 0\}} \frac{\Delta A(v_N,w_N)}{|||w_N|||}\\
&\leq \left(1+\frac{C_2}{C_1 C_5}\right) \widehat h^{m+t-s} \sum_{i=1}^{N} \|\wtu\|_{H^s(E_i)} + \frac{1}{C_1 C_5}\inf_{v_N\in V^N} \sup_{w_N\in V^N\backslash\{\bf 0\}} \frac{\Delta A(v_N,w_N)}{|||w_N|||},
\end{aligned}
\end{equation}
where we used \eqref{fe_convergence} and the fact that $\left.\wtu\right|_{E_i}=0$ for all $i=N+1,\ldots N+K$. 

Note that the first term of the right hand side of \eqref{HsStrang} depends on the grid size inside of $\omgs$; thus, it is not affected by a coarser nonuniform grid in $\omgc$. This is not the case for the second term which depends on the integration method and on the grid size in $\omgsc$. In the limit cases of very accurate integration formulas or of a fine grid in $\omgc$, the second term in \eqref{HsStrang} is negligible, whereas as the degree of accuracy of the numerical integration becomes smaller or as the grid becomes coarser, the leading term in the estimate is the ``error'' $\Delta A$. Another important factor to be considered is the fact that the integrand vanishes in regions of $\omgc$ very far from $\omgs$; thus, for such regions a coarser grid might not affect the accuracy of the solution. Hence, in order to find a bound for $\Delta A$ (and, as a consequence, a coarsening rule) one has to take into account the interaction of the accuracy of the quadrature rule, the grid size, and the structure of the integrand in $\omgc$. As of now, we do not have such a bound  nor an optimal coarsening algorithm, but we empirically choose larger $h_i$'s  far from $\omgs$. In \S\ref{oned}, we present an example of a coarsening algorithm which preserves accuracy, provided that a proper empirical tuning of its parameters is performed.

\section{Numerical tests}\label{sec:num_tests}

In this section we present the results of numerical tests for finite element discretizations of one-dimensional problems. Although preliminary in nature, these results allow us to illustrate the theoretical results presented in \S\S\ref{nonlocal_fractional} and \ref{fdapp} and to demonstrate the viability of using nonlocal diffusion problems as a vehicle for determining approximate solutions of fractional derivative problems posed on bounded domains. We first consider the convergence of approximate nonlocal solutions to the corresponding nonlocal solutions and we discuss the choice of a coarser grid in $\omgc$ in terms of accuracy and memory requirements. Then, we provide qualitative and quantitative comparisons between approximate nonlocal solutions and exact solutions of the fractional Laplacian equation.

\subsection{Finite element approximation of a one-dimensional problem}\label{oned}

We define a one-dimensional model problem and introduce its finite element approximation. Let $\omg=(-1,1)$ and $\omgc=(-1-\lambda,-1)\cup(1,1+\lambda)$ so that $\omgsc=(-1-\lambda,1+\lambda)$. In this case \eqref{wosd} becomes
\begin{equation}\label{wosd_1d}
\frac{c_{\ddd,s}}{2}\int_{-1-\lambda}^{1+\lambda}\int_{x-\lambda}^{x+\lambda} \frac{\widetilde u_N(y)-\widetilde u_N(x)}{|y-x|^{1+2s}} \big(v(y)-v(x)\big)\,dy\,dx= \int_{-1}^1 f v\,dx, \;\;\forall\,v\in V_c^N.
\end{equation}
For integers $K,N>0$, we introduce a partition of $\overline\omgsc = [-1-\lambda,  1+ \lambda]$ such that
\begin{equation}\label{eq:partition}
\begin{aligned}
&-1- \lambda= x_{-K}<\cdots<x_{-1}<0=x_0<x_1<\cdots<x_{N-1}\\&\qquad\qquad
<x_{N}=1<x_{N+1}<\cdots.<x_{N+K} = 1+\lambda.
\end{aligned}
\end{equation} 
Then, we define $h$ as $\max_{i = -K,\ldots,K+N-1}|x_{i+1}-x_i|$. We choose $V_c^N(\omgsc)$ to be the finite element space of continuous piecewise-linear polynomials with respect to the partition \eqref{eq:partition}. As a basis, we choose the standard ``hat'' functions $\{\phi_j(x)\}_{j=-K}^{K+N}$. Let $\wtuN(x)=\sum_{j=-K}^{K+N}U_j\phi_j(x)$ for some coefficients $U_j$; because the hat functions satisfy $\phi_{i}(x_j)=\delta_{ij}$, we have that $U_j = \widetilde u_N(x_j)$. Then, \eqref{wosd_1d} is equivalent to
\begin{equation}\label{weakFE1d}
\begin{aligned}
\frac{c_{\ddd,s}}{2} \sum_{j=1}^{N-1}
 U_j \hspace{-2pt}  \int_{-1-\lambda}^{1+\lambda}
 \hspace{-3pt}\int_{x-\lambda}^{x+\lambda} \frac{\phi_j(y)-\phi_j(x)}{|y-x|^{1+2s}} \big(\phi_{i}(y)-\phi_{i}(x)\big)\,dydx = 
\int_{-1}^1 f\phi_i \,dx.,\\
i=1,\ldots,N-1,
\end{aligned}
\end{equation}
Also, because $\widetilde u=0\in \omgc$, $U_0=U_N=0$.

The assembly of the matrix and right-hand side of the linear system corresponding to \eqref{weakFE1d} requires the use of quadrature rules for the approximate evaluation of the integrals involved. We first break up all integrals into sums of integrals over the subintervals $(x_i,x_{i+1})$ of the partition \eqref{eq:partition}. For the single integrals over $\omgs=(-1,1)$ in \eqref{weakFE1d}, we use a four-point Gauss quadrature rule with respect to each subinterval. The double integrals must be treated with greater care, due to the singularity of the integrand. For the inner integrals, we split the integral over $(x- \lambda, \, x+\lambda)$ into the sum of integrals over $(x- \lambda, \, x)$ and $(x, \, x+\lambda)$, then split each of those integrals into a sum over the subintervals of the grid, and then use a four-point Gauss quadrature rule over each subinterval; in this way, we avoid the singularity at $x$. An accurate and efficient choice of quadrature rule for the outer integral is a more delicate issue. For this purpose, we use the built-in Matlab function {\tt quadgk.m} which implements a high-order global adaptive Gauss-Kronrod quadrature formula which is well-suited for integrands having mild end-point singularities \cite{Shampine}. Another advantage of using this Matlab function is that all subintervals are processed at the same time so that the number of function evaluations is reduced. 

Note that the double integrals on the left-hand side correspond to a nonlocal finite element stiffness matrix; an important issue in the numerical solution of nonlocal problems is related to its structure. In fact, differently from the local case where the matrix is sparse and banded, in the nonlocal case the matrix is dense and the inner integral over $(x- \lambda, \, x+\lambda)$ spans a large part of $\omgsc$, especially when $\lambda\gg I$, which is often the case. Thus, computations are much more onerous compared to that for differential equations. 

\subsection{The analytic solution of the fractional Laplacian}
\begin{figure}[h]
\begin{center}
\begin{tabular}{cc}
\includegraphics[width=0.4\textwidth]{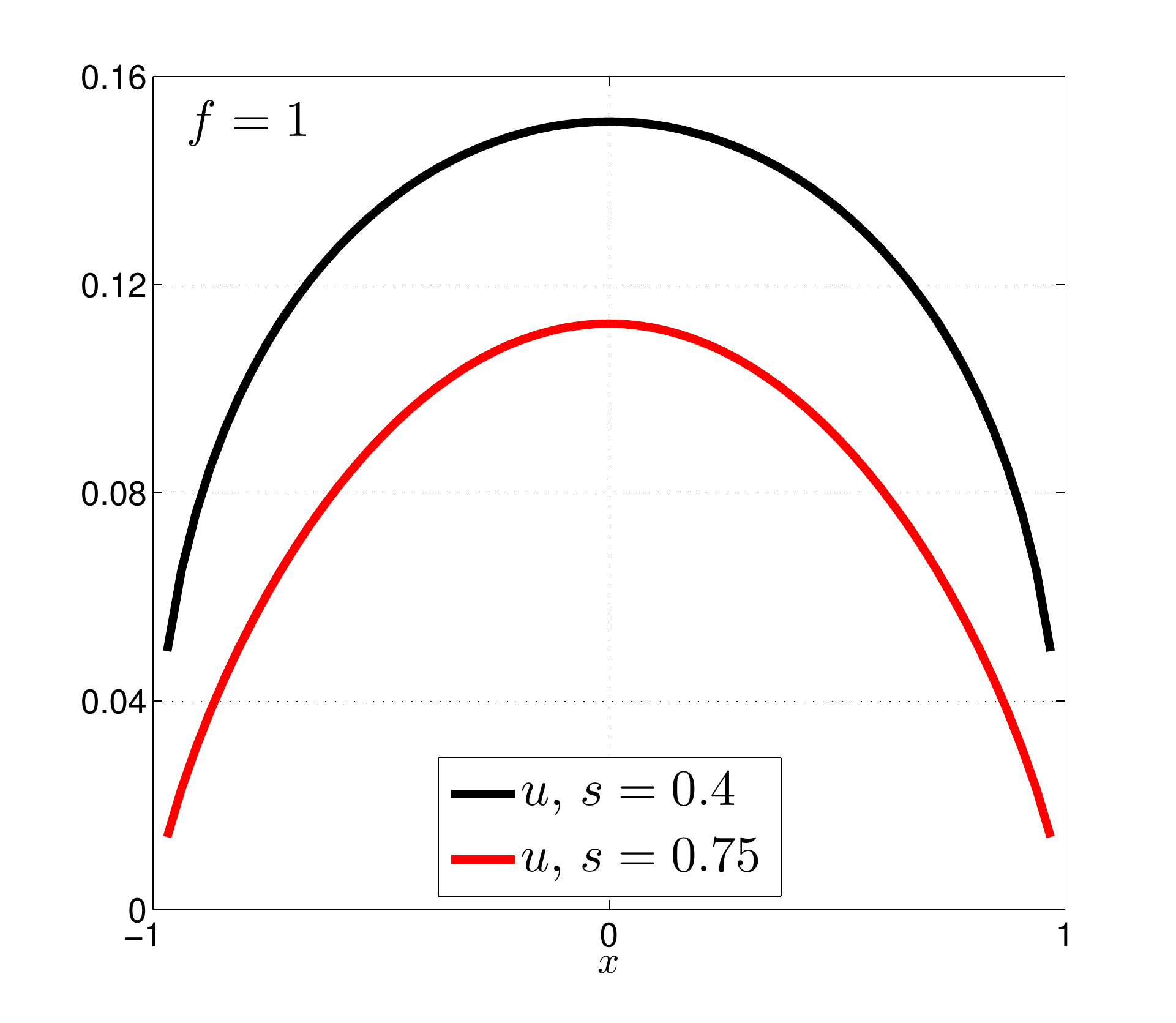} &
\includegraphics[width=0.4\textwidth]{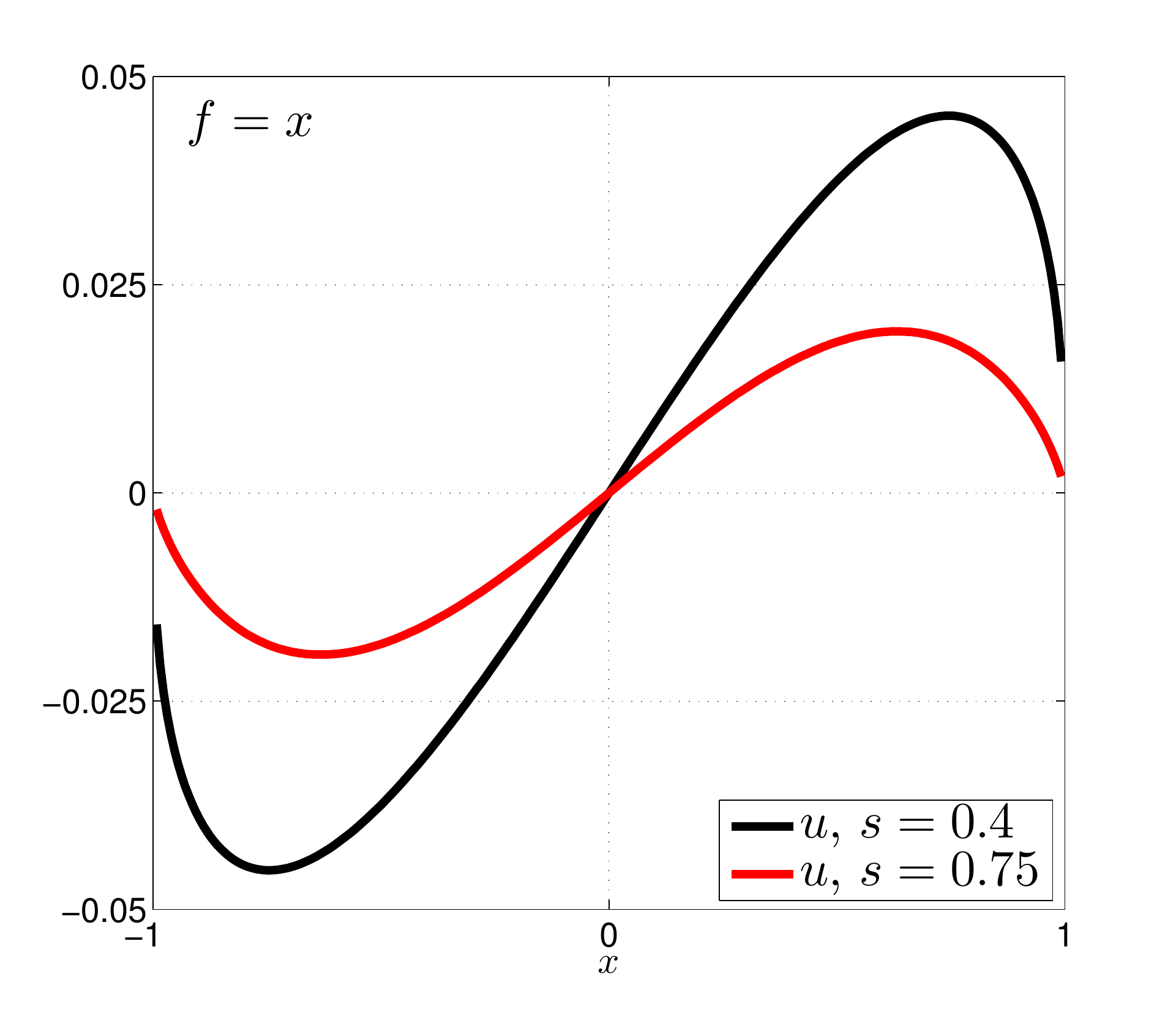}
\end{tabular}
\end{center}
\caption{Analytic solution $u$ for cases Ia,b in \eqref{data1} (left) and cases IIa,b in \eqref{data2} (right).}
\label{fig:anl_sol}
\end{figure}

We consider the analytic solution of problem \eqref{eq:fl_strong} in the ball $B_R(\bf 0)$ of radius $R$ in $\mbRn$. For every function $f\in L^2(B_R(\bf 0))$ there exists a unique function $u\in H^s(\mbRn)$, given explicitly using the Green's function \cite{diblasio} by
\begin{equation}\label{eq:int_u}
u(\xb) = \int_{B_R(\bf 0)} G(\xb,\yb) f(\yb)\,d\yb,
\end{equation}  
where 
\begin{equation}\label{eq:greens}
G(\xb,\yb)=C_{\ddd,s} |\xb-\yb|^{2s-n}\int_0^{r_0(\xb,\yb)} \dfrac{r^{s-1}}{(r+1)^{n/2}}dr,
\end{equation}
with 
\begin{equation}\label{eq:Cns_r0}
C_{\ddd,s}=s^{-2s}\dfrac{\Gamma\left(\frac{n}{2}\right)}{R^2\pi^{n/2}\Gamma(s)^2} 
\qquad {\rm and} \qquad 
r_0(\xb,\yb)=\dfrac{(R^2-|\xb|^2)(R^2-|\yb|^2)}{|\xb-\yb|^2}.
\end{equation}
In the particular case of $f=1$ in $B_1(\xb)$ 
\begin{equation}\label{eq:intu1}
u(\xb)=2^{-2s}\,\dfrac{\Gamma\left(\frac{n}{2}\right)}{\Gamma\left(\frac{n+2s}{2}\right)\Gamma(1+s)}\, (1-|\xb|^2)^s.
\end{equation}
In the one-dimensional numerical tests we consider the following data sets:
\begin{subequations}\label{cases1}
\begin{align}
\mbox{Ia }\,\,& 
 \left\{\begin{aligned}	
	& f(x) = 1\\
	& s = 0.75
\end{aligned}\right. &
\mbox{Ib }\,\,& 
 \left\{\begin{aligned}	
	&  f(x) = 1\\
	& s = 0.40
\end{aligned}\right.\label{data1}
\\[1ex]
\mbox{IIa }\,\,& 
 \left\{\begin{aligned}	
	& f(x) = x\\
	& s = 0.75
\end{aligned}\right. &
\mbox{IIb }\,\,& 
 \left\{\begin{aligned}	
	& f(x) = x\\
	& s = 0.40.
\end{aligned}\right.\label{data2}
\end{align}
\end{subequations}
Note that we choose values of $s$ both less than and greater than $1/2$. The corresponding analytic solutions are shown in Figure \ref{fig:anl_sol}.

\subsection{Convergence of finite element approximations to the nonlocal solution}
We examine the convergence, with respect to the grid size $h$, of finite element approximations. These results are not related to the solution of the fractional Laplacian equation and they only serve to illustrate the theoretical results about the convergence of the approximate nonlocal solutions $\wtuN$ to the solution $\wtu$ of \eqref{eq:eps_weak}. We consider the data sets Ia and Ib; because an exact solution is not available, we define a surrogate $\wtu_R$ for $\wtu$ to be the finite element solution on a very fine uniform grid. We then report on the ``error'' $\|\wtu_R-\wtuN\|_{L^2(\omgs)}$, i.e., the difference between the surrogate and the approximation $\wtuN$ obtained using coarser uniform grids and the same interaction radius $\lambda$. In particular, the surrogate is determined using\footnote{This small interaction radius has been chosen reduce the computational costs compared to that which would be incurred for larger values of $\lambda$; in fact, in general, with uniform grids on $\omgsc$, one can afford to use only small interaction radii.} $\lambda=0.1$ and $h=2^{-11}$. In Table \eqref{tab:hrate}, we list the error and the corresponding rate. In both cases we observe a (optimal) rate of of convergence of approximately $0.5$. This is expected; the analytic solution belongs to $H^{0.5-\varepsilon}(\mbR)$; in fact, it is continuous but its derivative has an infinite discontinuity at $x=-1,1$. As a consequence, the estimate \eqref{fe_convergence} does not hold because $u\notin H^{m+t}(\mbR)$ with $m=1$ and $t\in[s,1]$ and we cannot expect a rate better than 0.5.

\begin{table}[h!]
\begin{center}
\begin{tabular}{ c | c | c | c | c | }
\cline{2-5}
								&	\multicolumn{2}{|c|}{Ia}		& 	\multicolumn{2}{|c|}{Ib}\\
\hline                  
\multicolumn{1}{|c|}{$h$}		&	error		&	rate 		&	error			&	rate\\	
\hline
\multicolumn{1}{|c|}{$2^{-3}$}	&	6.92e-02		& 	- 			&	6.23e-02				& 	-	\\
\hline
\multicolumn{1}{|c|}{$2^{-4}$}	&	4.74e-02 	& 	0.55			&	4.53e-02				& 	0.46	\\
\hline
\multicolumn{1}{|c|}{$2^{-5}$}	&	3.17e-02 	&	0.58			&	3.08e-02  			& 	0.55	\\
\hline
\multicolumn{1}{|c|}{$2^{-6}$}	&   2.11e-02 	&	0.58			&	2.08e-02 			& 	0.55	\\
\hline
\end{tabular}
\caption{Dependence on the grid size $h$ of the ``error'' $\|\wtu_R-\wtuN\|_{L^2(\omgs)}$ and the rate of convergence of approximate nonlocal solutions. Results are provided for the test cases Ia and Ib.}
\label{tab:hrate}
\end{center}
\end{table}

\paragraph{Nonuniform grid in $\omgc$} For large values of $\lambda$, using a uniform grid in $\omgc$ is not computationally affordable so that using a coarser grid is mandatory. However, this might affect the accuracy of the solution; thus, we need a coarsening algorithm that preserves accuracy. Of course, this is not a trivial task (see \S\ref{sec:non_unif}). We design an algorithm which is not optimal but still allows us to maintain a certain level of accuracy with a proper tuning of parameters. 

\begin{figure}[h]
\begin{center}
\begin{tabular}{c}
\includegraphics[width=0.80\textwidth]{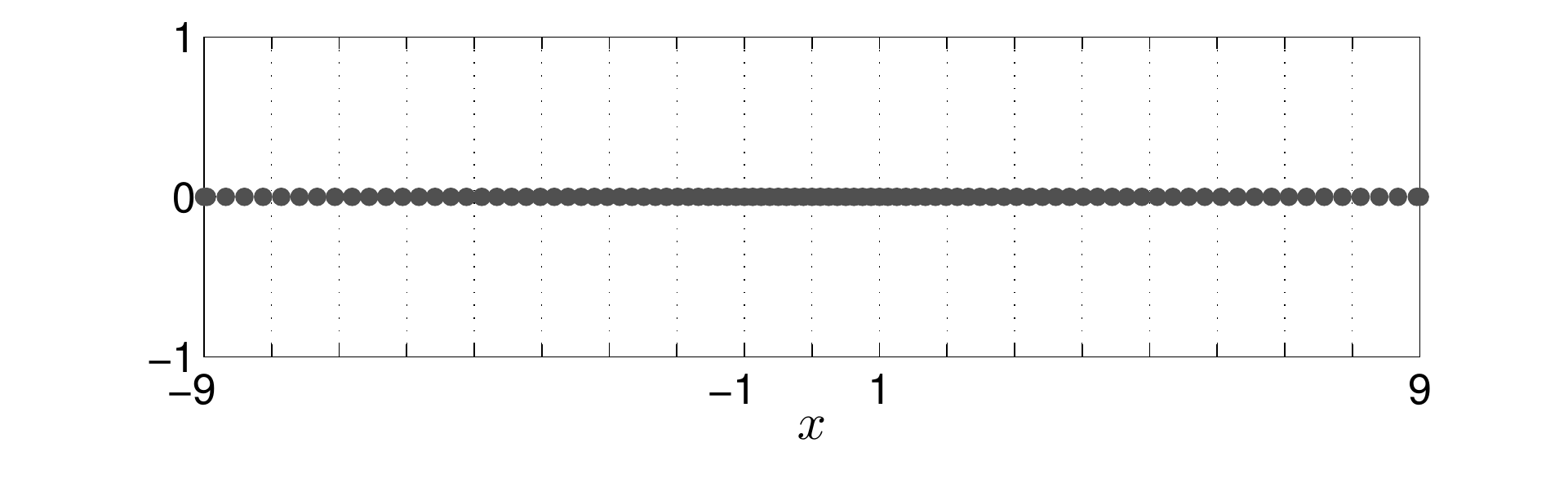}\\
\includegraphics[width=0.80\textwidth]{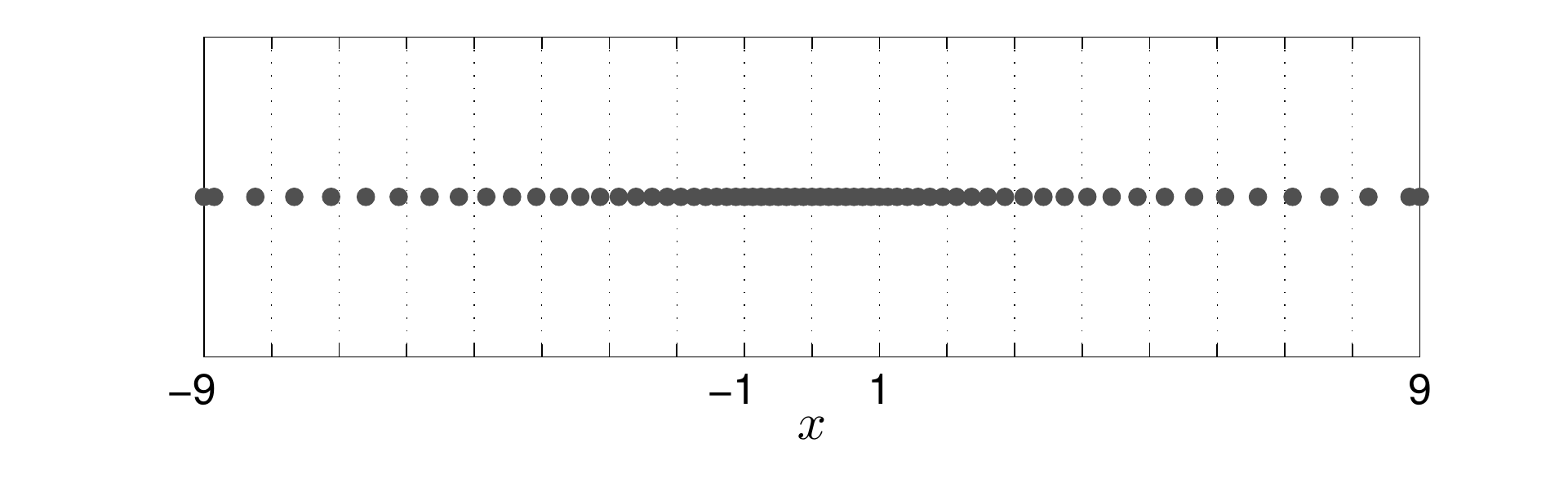}\\
\includegraphics[width=0.80\textwidth]{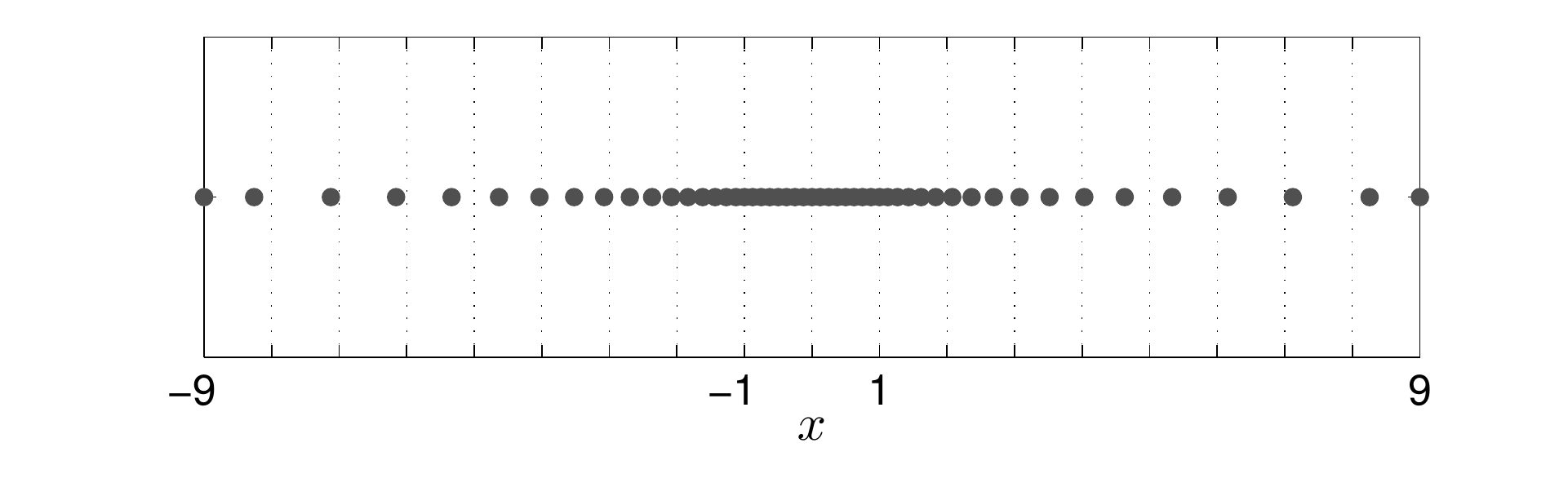}
\end{tabular}
\end{center}
\caption{Three grid configurations for $\lambda=2^3$, $N=2^4$, and $p=0.5,1,1.5$, respectively, from top to bottom.}
\label{fig:nodes}
\end{figure}

Referring to the partition \eqref{eq:partition}, we let $\widehat{x}=(x_0+x_N)/2$ and for $i=-K,\ldots -1$ and $i=N+1,\ldots N+K$, we define $x_i$ as 
\begin{equation}\label{eq:coarse_partition}
x_i = \left\{
\begin{array}{ll}
x_{i+1} - \widehat{h}\left(\dfrac{\widehat{x}-x_{i+1}}{\widehat{x}-x_0}\right)^p & i=-1,\ldots -K\\
[3.5mm]
x_{i-1} + \widehat{h}\left(\dfrac{x_{i-1}-\widehat{x}}{x_N-\widehat{x}}\right)^p & i=N+1,\ldots N+K,
\end{array}\right.
\end{equation}
where $K$ is chosen so that $x_{-K+1}>-1-\lambda$ and $x_{N+K-1}<1+\lambda$; then, $x_{-K} = -1-\lambda$ and $x_{N+K}=1+\lambda$; $\widehat{h}$ is the grid size inside of $\omgs$ and $p$ is the parameter that determines the coarseness: the larger $p$, the coarser the grid. For some values of $N$ and $\lambda$, we test the convergence of the solutions obtained with a nonuniform grid to that obtained with a uniform one for several values of $p$. In Table \ref{tab:coarsening}, we report the dimension of the partition for $N=2^7$ and $\lambda=2^2,2^3$; $p=0$ corresponds to the uniform grid and it is reported as a reference value. Unfortunately, the difference between the solutions on nonuniform and uniform grids, though decreasing for decreasing values of $p$, does not have a specific convergence rate. For the numerical tests presented in the following sections, $p$ is chosen empirically in such a way that the coarsening does not affect the accuracy, i.e., $|\Delta A(w_n,v_n)|/|||w_N|||$ in \eqref{HsStrang} is negligible for all $w_N,v_N\in V^N$. Figure \ref{fig:nodes} shows, for $\lambda=2^3$ and $N=2^4$, the computational grids for $p=0.5,1,1.5$.

\begin{table}[h!]
\begin{center}
\begin{tabular}{ | c | c | c |  }
\hline                  
$p$		&	$\lambda=2^2$	&	$\lambda=2^4$	\\	\hline
1.000	&	337				&	413				\\	\hline
0.500	&	447				&	643				\\	\hline
0.250	&	531				&	847				\\	\hline
0.125	&	583				&	985				\\	\hline
0.000	&	641				&	1153				\\	\hline
\end{tabular}
\caption{Dependence on the parameter $p$ of the dimension of the grid for $N=2^7$. In this case, for $p<0.5$ the coarsening error is negligible.}
\label{tab:coarsening}
\end{center}
\end{table}

\subsection{Comparison of $u$ and $\widetilde u$}
\begin{figure}[h]
\begin{center}
\begin{tabular}{cc}
\includegraphics[width=0.48\textwidth]{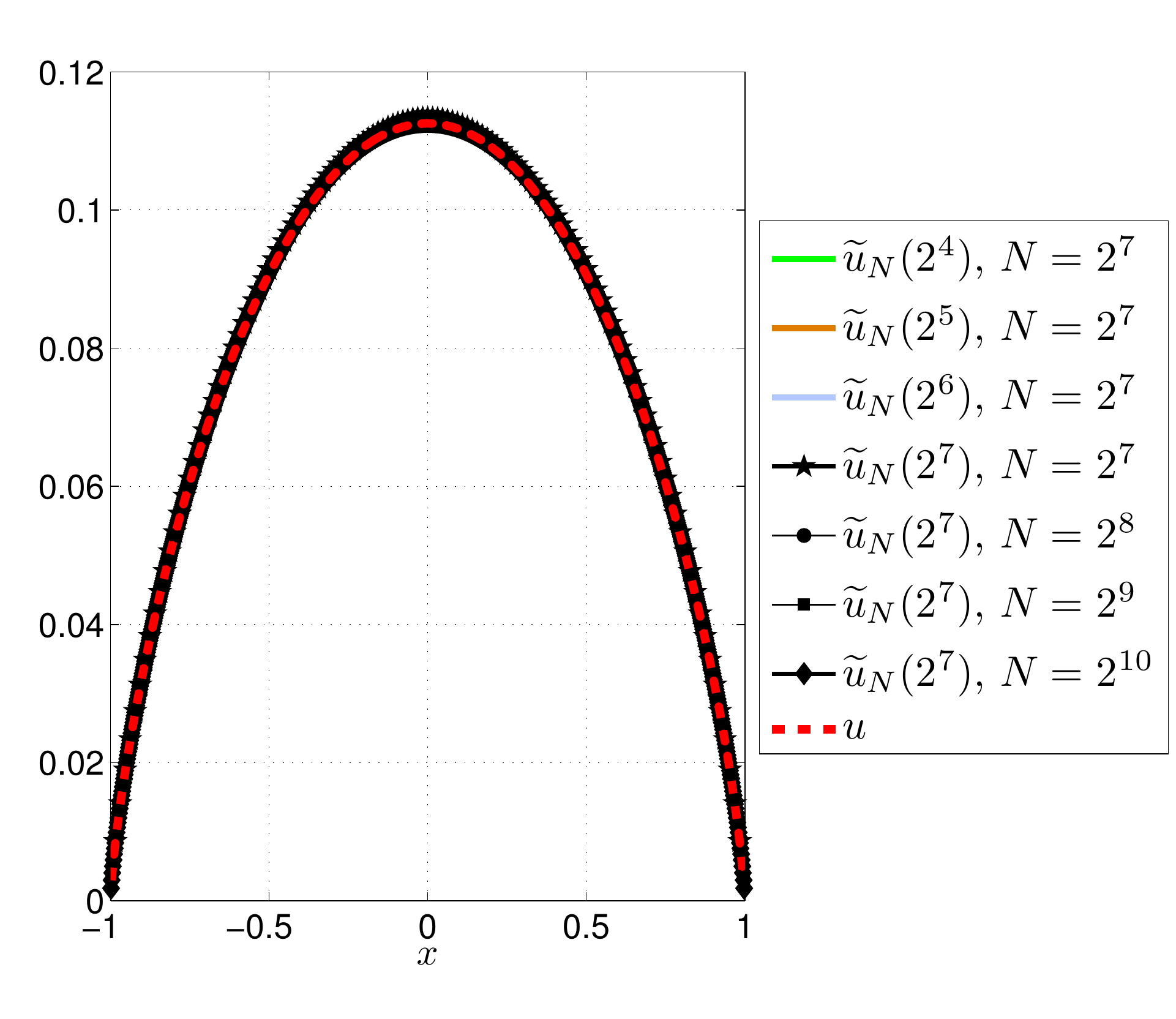} &
\includegraphics[width=0.45\textwidth]{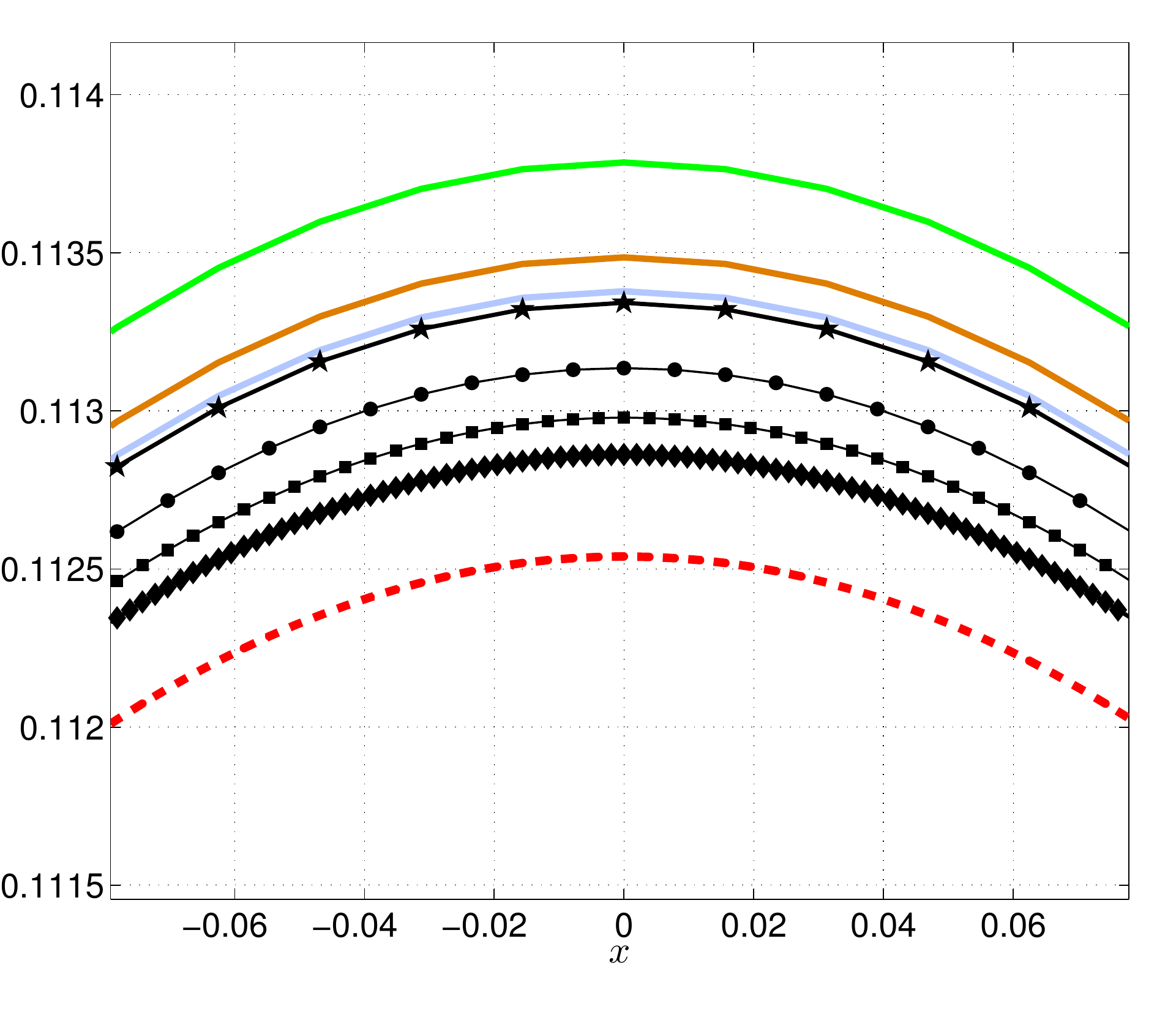}\\
\includegraphics[width=0.48\textwidth]{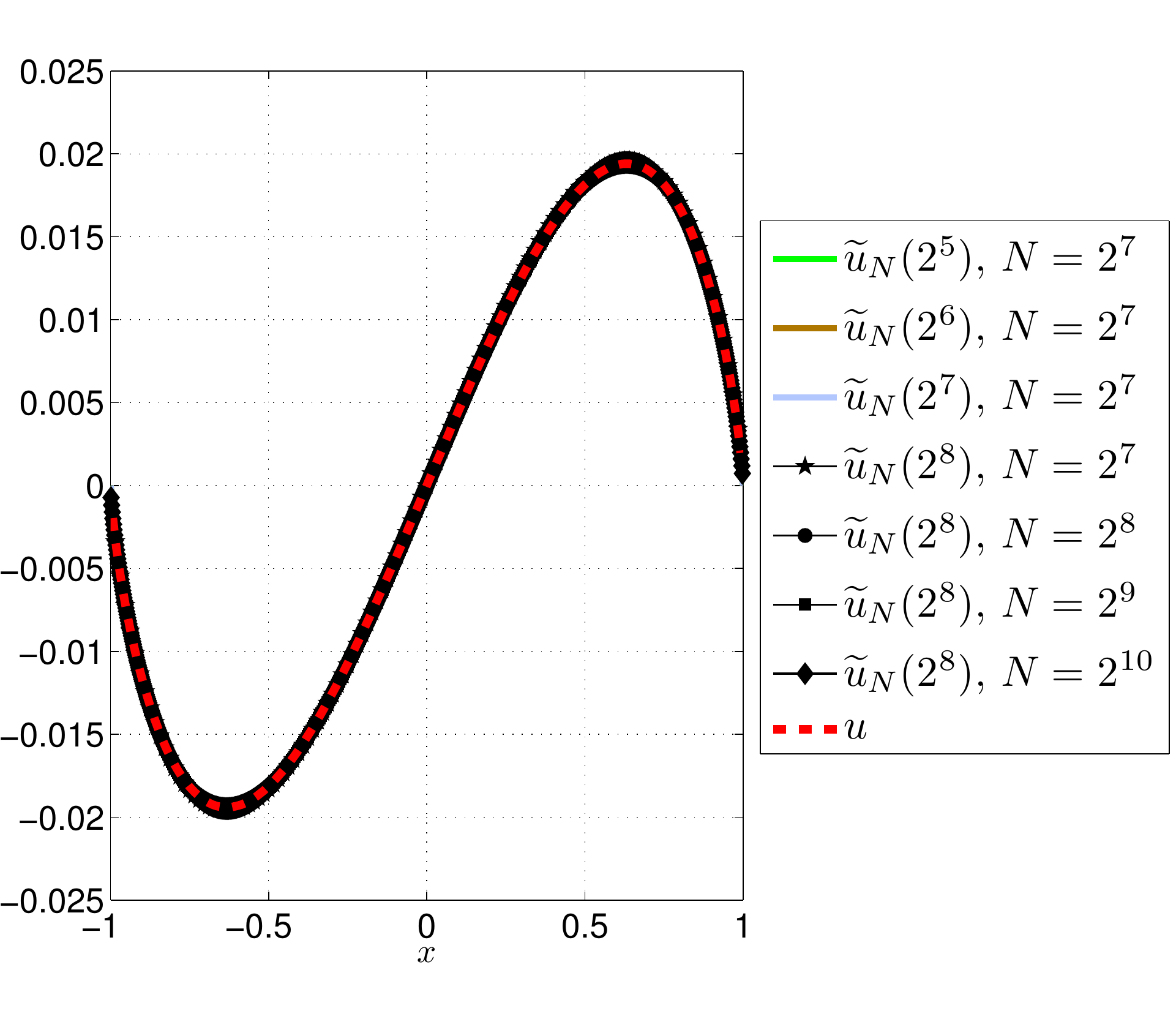} &
\includegraphics[width=0.45\textwidth]{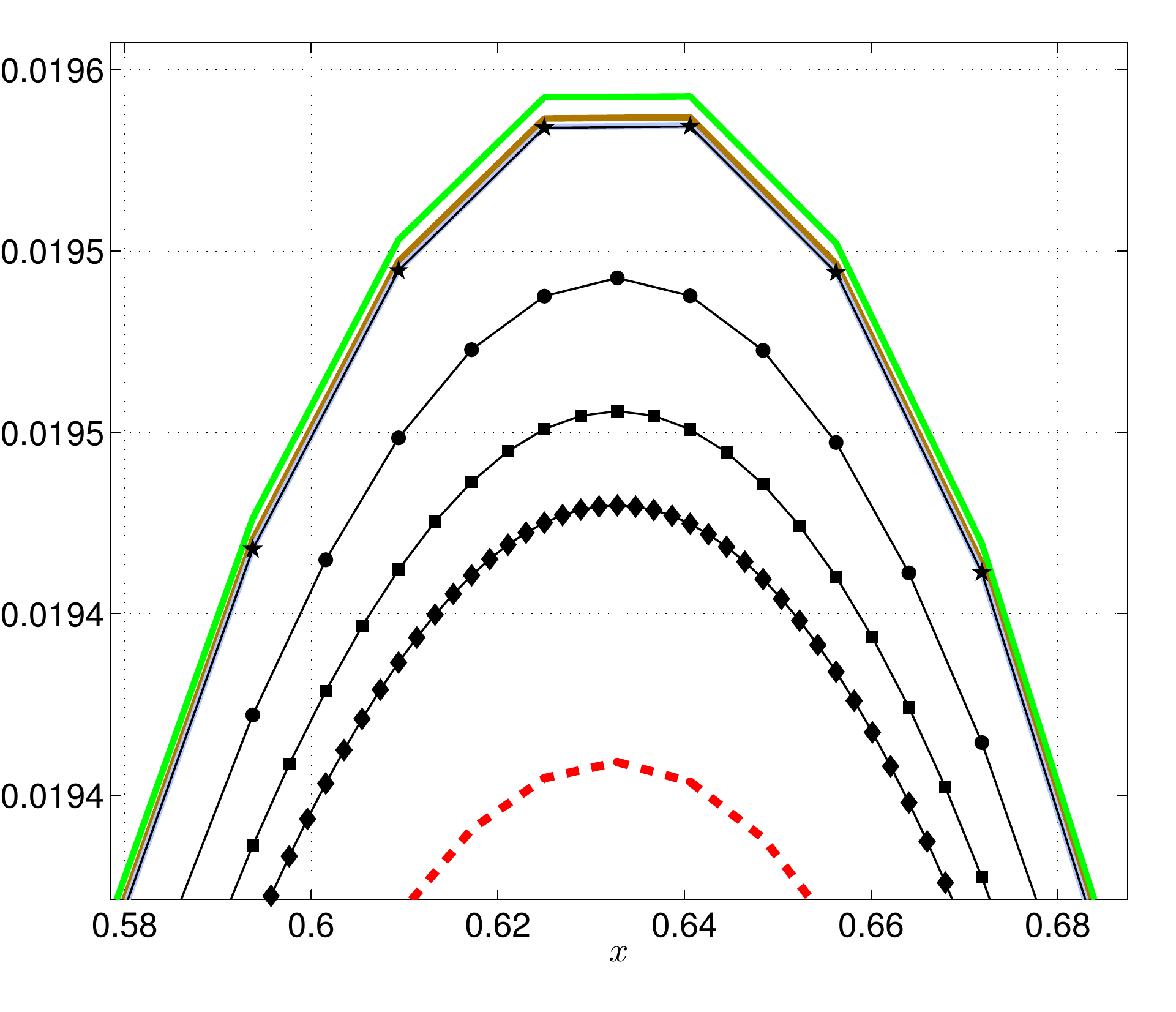}
\end{tabular}
\end{center}
\caption{Numerical solutions $\widetilde u_{N,\lambda}$ for different values of $N$ and $\lambda$ show the convergence to the solution $u$ of the fractional Laplacian as $N\to \infty$ and $\lambda\to\infty$. On the right, zoom-ins of the solutions peaks. Results are provided for the test cases Ia (top) and IIa (bottom).}
\label{fig:loc_lim75}
\end{figure}

Through several graphical examples we show how accurately the nonlocal solution can approximate the solution of the fractional Laplacian equation posed on bounded domains. In Figure \ref{fig:loc_lim75}, left column, we report the finite element solutions $\widetilde u_{N,\lambda}$ and the analytic solution $u$ for the data sets Ia and IIa using several values of the interaction radius $\lambda$ and of the grid dimension $N$. To appreciate the differences of the solutions we report, in the right column, zoom-ins near the peaks of the analytic solutions. For a fixed value of $N$ and increasing values of $\lambda$, the approximate nonlocal solutions converge to $u$ as predicted by theory; we also observe a convergence to $u$ for increasing values of $N$. For all $s>0.5$, similar results are obtained, i.e., the larger $\lambda$ and $N$, the closer is $\widetilde u_{N,\lambda}$ to $u$. 

\begin{figure}[h]
\begin{center}
\begin{tabular}{c}
\includegraphics[width=0.9\textwidth]{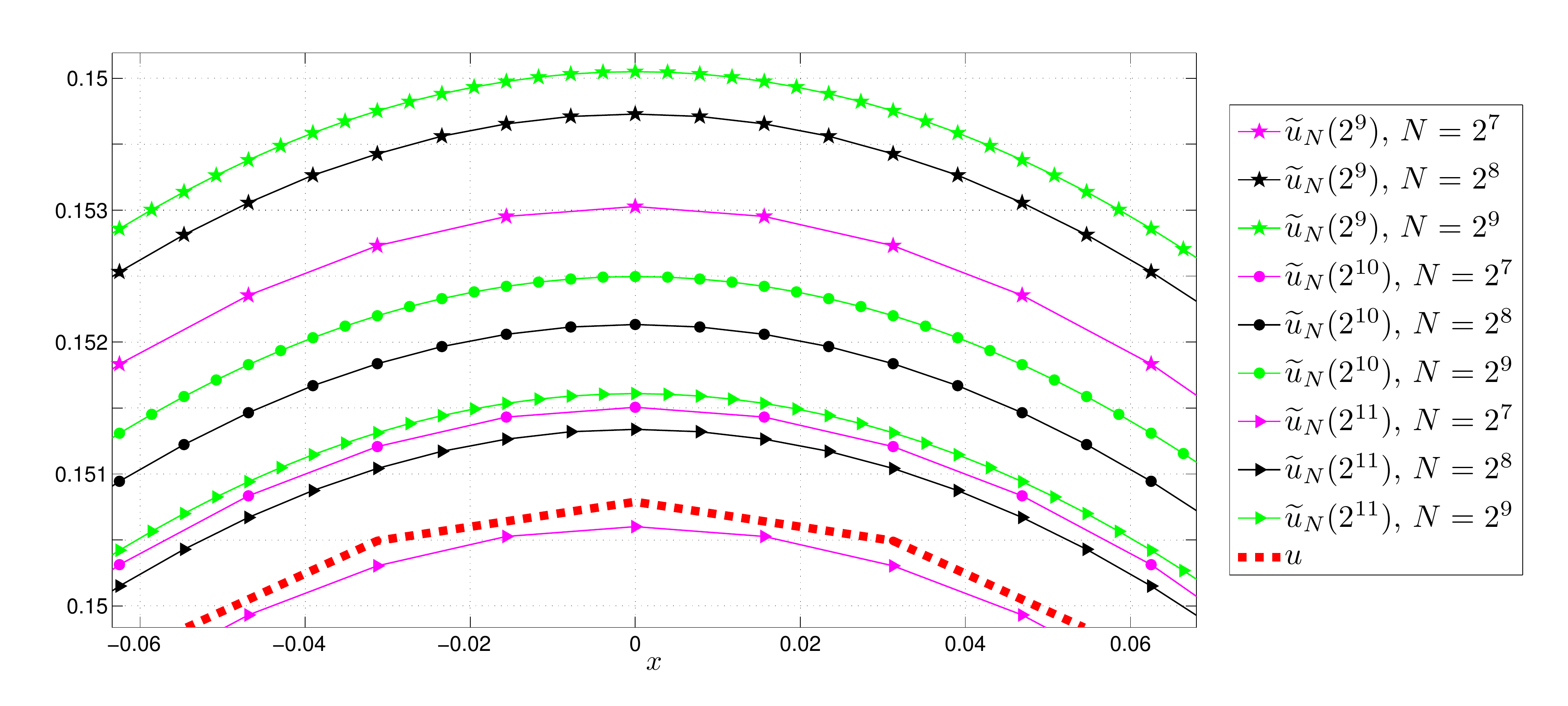} \\
\includegraphics[width=0.90\textwidth]{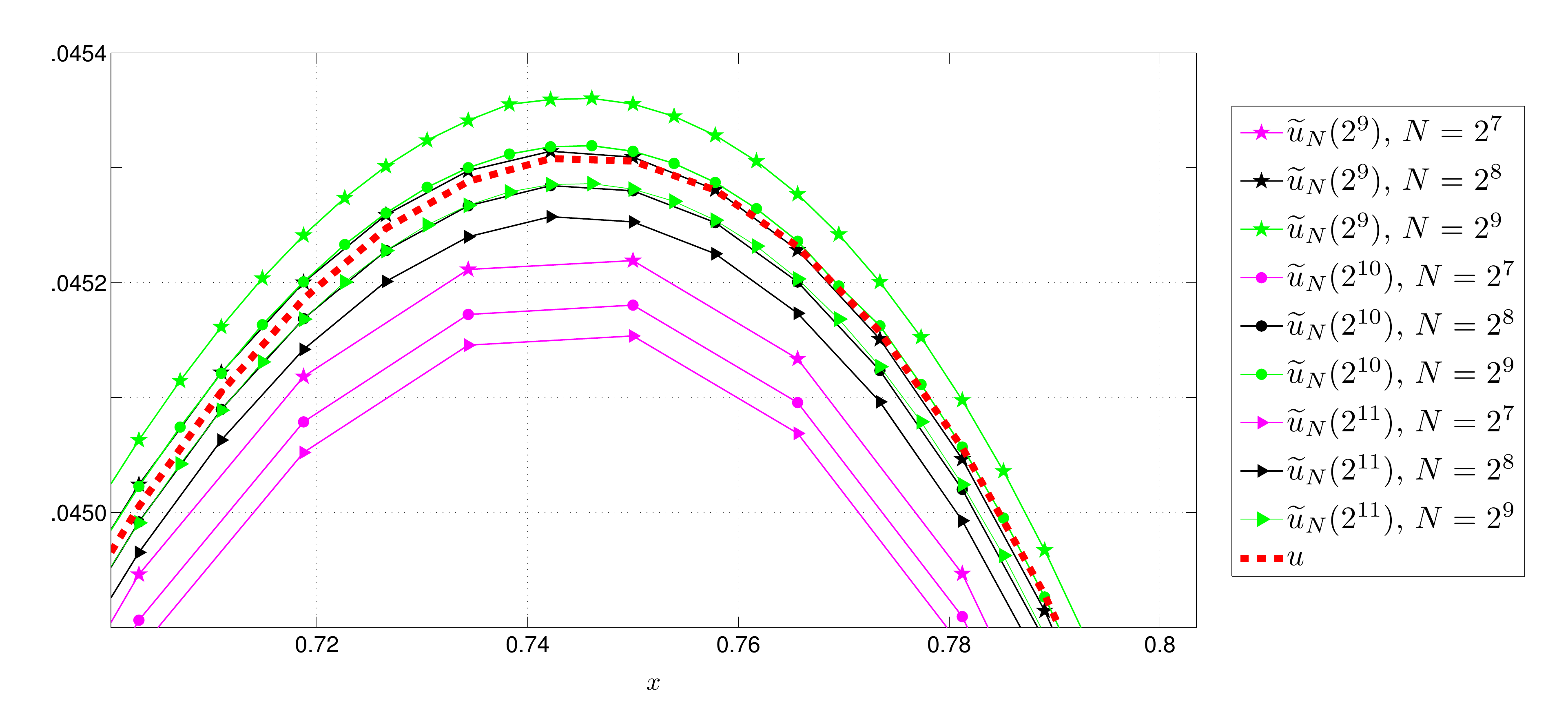}
\end{tabular}
\end{center}
\caption{Numerical solutions $\widetilde u_{N,\lambda}$ for different values of $N$ and $\lambda$. The figures are zoom-ins of the solutions peaks. Results are provided for the test cases Ib (top) and IIb (bottom).}
\label{fig:loc_lim40}
\end{figure}

\begin{figure}[h]
\begin{center}
\begin{tabular}{c}
\includegraphics[width=0.9\textwidth]{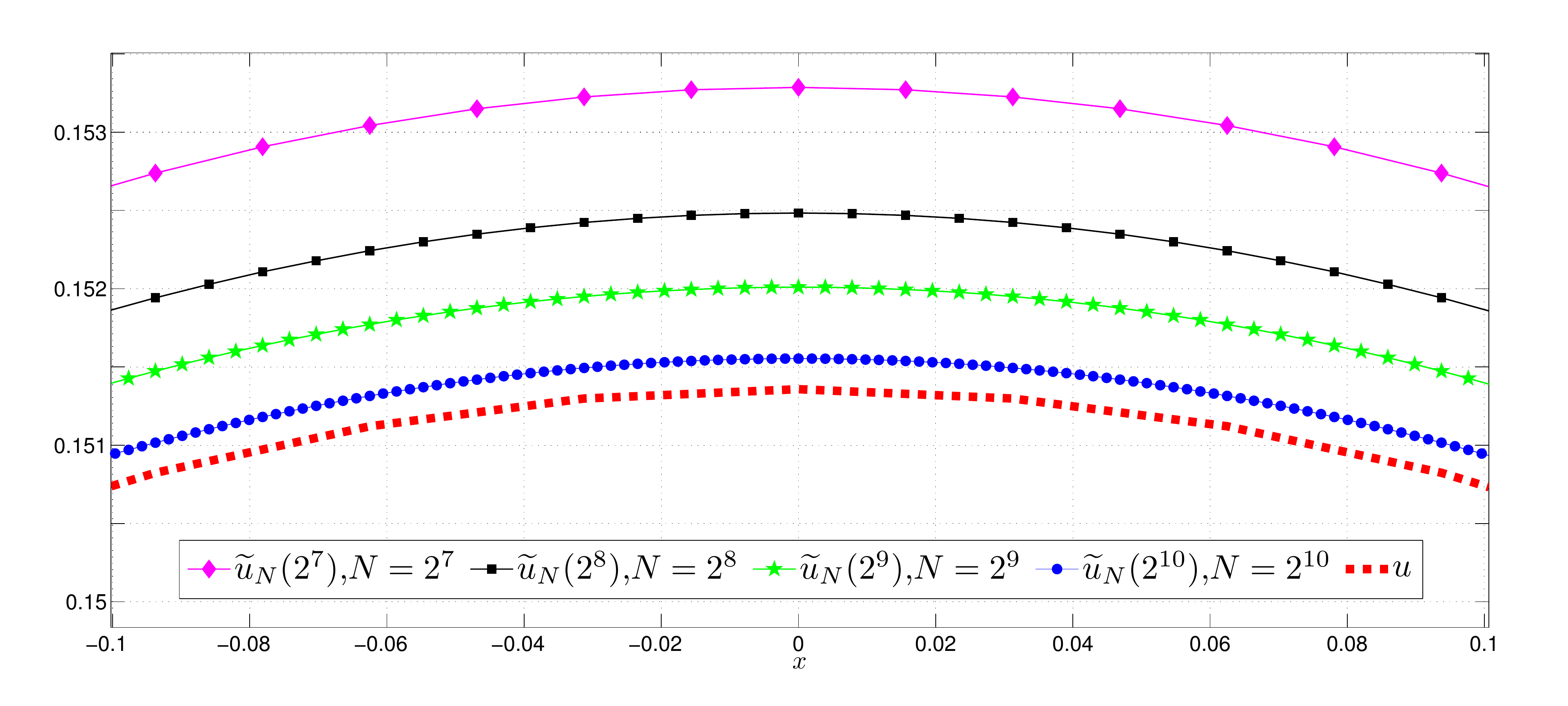} \\
\includegraphics[width=0.9\textwidth]{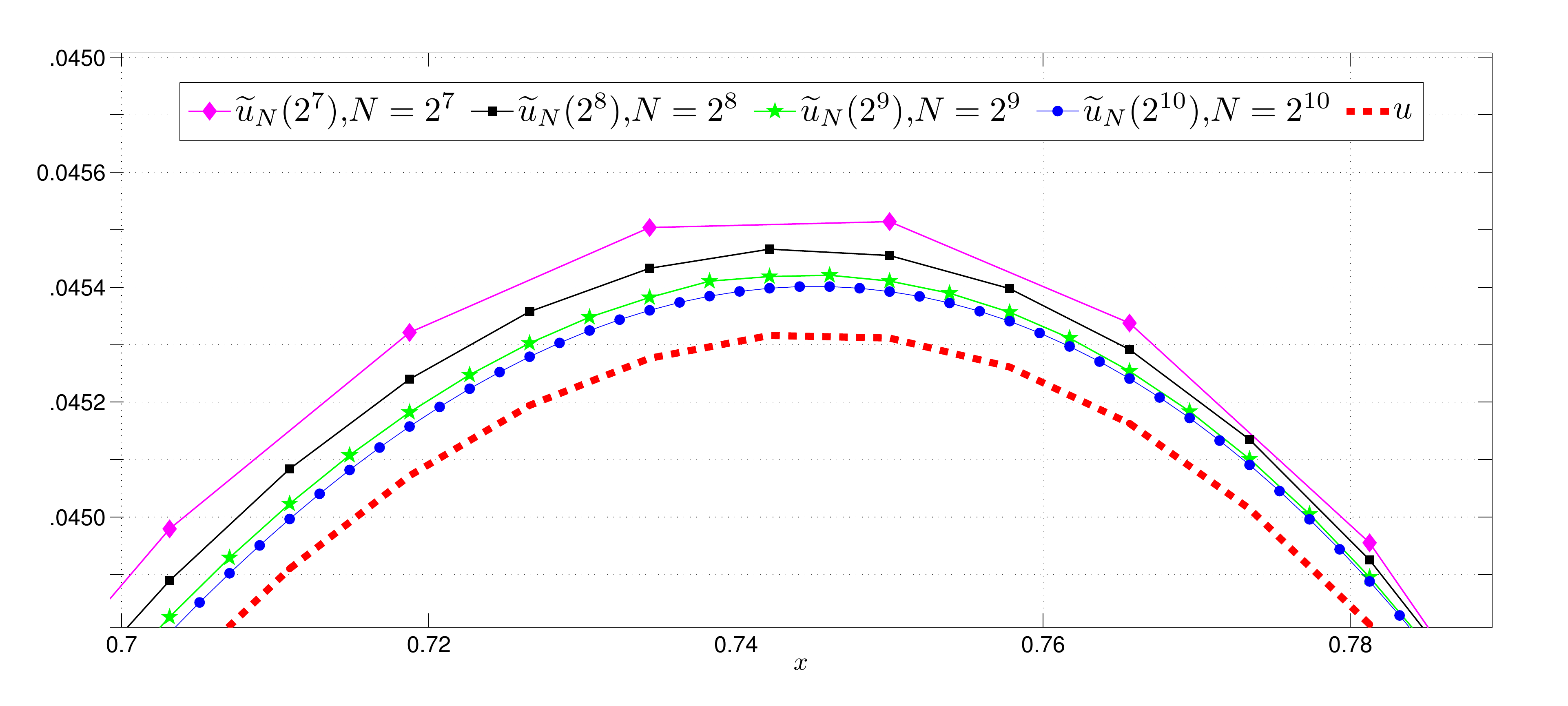}
\end{tabular}
\end{center}
\caption{Numerical solutions $\widetilde u_N$ for different values of $h$ and $\lambda$ show the convergence to the solution $u$ of the fractional Laplacian as $N\to \infty$ and $\lambda\to\infty$. The figures are zoom-ins of the solutions peaks. Results are provided for the test cases Ib (top) and IIb (bottom).}
\label{fig:loc_lim40DC}
\end{figure}

For $s<0.5$, the situation is quite different and deserves a more detailed analysis. In Figure \ref{fig:loc_lim40}, we report $\widetilde u_{N,\lambda}$ and $u$ for $\lambda=2^9,2^{10},2^{11}$ and $N=2^7,2^8,2^9$; again, the plots are zoom-ins around the peaks of $u$ for the data sets Ib (top) and IIb (bottom). Here, changes in $\lambda$ and $N$ have different effects: increasing values of $\lambda$ lead to an amplitude reduction whereas increasing values of $N$ lead to an amplitude increase. We conjecture that for a fixed very small $h$, it is possible to observe a convergence to $u$ from above for increasing values of $\lambda$. Because of computational restrictions we cannot utilize a grid which is fine enough to observe such behavior. For this reason we find it significant to study the convergence of $\widetilde u_{N,\lambda}$ to $u$ for increasing values of $\lambda$ and $N$ simultaneously; in Figure \ref{fig:loc_lim40DC}, we then observe a convergence to $u$ from above.

\subsection{Rate of convergence to the solution of the fractional Laplacian equation}\label{sec55}
We study the convergence of the approximate nonlocal solution $\widetilde u_{N,\lambda}$ to the analytic solution $u$ of the fractional Laplacian equation.  Recall that the error bound depends on both $\lambda$ and $N$ (see \eqref{u_uNlambda_estimate}); through several examples, we show how the choice of $\lambda$ and $N$ affects the accuracy of the numerical solution and, as a consequence, how it affects the rate of convergence. First, we consider simultaneously increasing $\lambda$ and $N$; in Table \ref{tab:IAIBrate}, we report the approximation error, i.e., $\|u-\widetilde u_{N,\lambda}\|_{L^2(\omgs)}$, and its rate. In this experiment $p=0.5$, except for $h=2^{-8}$ where $p=0.6$. We observe a rate of $0.5$, the same as the one obtained for the finite element approximations of nonlocal solutions. We conclude that the convergence with respect to $u$ is dominated by the finite element approximation error in \eqref{u_uNlambda_estimate}.
 
\begin{table}[h!]
\begin{center}
\begin{tabular}{ c  c | c | c | c | c |  }
\cline{3-6}
							&								&\multicolumn{2}{|c|}{Ia}	& 	\multicolumn{2}{|c|}{Ib}\\	\hline                  
\multicolumn{1}{|c|}{$h$}	&\multicolumn{1}{c|}{$\lambda$}	&	error		&	rate 	&	error		&	rate		\\	\hline
\multicolumn{1}{|c|}{$2^{-4}$}	&\multicolumn{1}{c|}{$2^3$}		&	2.82e-02 	& 	- 		&	4.47e-02		& 	-		\\	\hline
\multicolumn{1}{|c|}{$2^{-5}$}	&\multicolumn{1}{c|}{$2^4$}		&	2.00e-03 	& 	0.50		&	2.98e-02		& 	0.59		\\	\hline
\multicolumn{1}{|c|}{$2^{-6}$}	&\multicolumn{1}{c|}{$2^5$}		&	1.41e-02 	&	0.50		&	2.03e-02  	& 	0.55		\\	\hline
\multicolumn{1}{|c|}{$2^{-7}$}	&\multicolumn{1}{c|}{$2^6$}		&   9.97e-03 	&	0.50		&	1.40e-02 	& 	0.53		\\	\hline
\multicolumn{1}{|c|}{$2^{-8}$}	&\multicolumn{1}{c|}{$2^7$}		&   7.05e-03 	&	0.50		&	9.76e-03 	& 	0.52		\\	\hline
\end{tabular}
\caption{Dependence on the grid size $h$ and on the interaction radius $\lambda$ of the error $\|u-\wtu_{N,\lambda}\|_{L^2(\omgs)}$ and the rate of convergence of approximate nonlocal solutions. Results are provided for the test cases Ia and Ib.}
\label{tab:IAIBrate}
\end{center}
\end{table}
In order to observe the convergence with respect to $\lambda$, we consider a fixed grid and we compare $\widetilde u_{N,\lambda}$ with $u$; as we conjectured in the previous section, though we do observe a convergence, a (not affordable) very fine grid is required to observe the $-2s$ rate. For this reason, we do not find it significant to report those results. We examine, instead, the convergence, for a fixed $N$, of the nonlocal approximate solution $\widetilde u_{N,\lambda}$ to a surrogate $u_R$ defined by solving for a finite element approximation using a very fine grid and a very large $\lambda$. In Tables \ref{tab:Arate} and \ref{tab:Brate} we report, for the test cases I and II, the ``errors'' $|||e|||=|||u_R-\widetilde u_{N,\lambda}|||$ and $\|e\|_{L^2(\omg)}=\|u_R-\widetilde u_{N,\lambda}\|_{L^2(\omgs)}$, i.e., the energy and $L^2$ norms of the difference between the surrogate and the finite element approximation using the same grid and smaller radius. Specifically, the surrogate $u_R$ for the analytic solution is determined using $h=2^{-9}$, $\lambda = 2^{11}$, and $p=1.25$. We observe the rate $-2s$ predicted by the estimate \eqref{eq:lambda_limit}, for both the energy norm and the $L^2$ norm. Note that because the observed rates of convergence with respect to the two norms are the same, the sharpness of the $L^2$ error estimate in \eqref{eq:lambda_limit} is seemingly sharp.

\begin{table}[h!]
\begin{center}
\begin{tabular}{ c | c | c | c | c | c | c | c | c | }
\cline{2-9}
								&	\multicolumn{2}{|c|}{Ia}		& 	\multicolumn{2}{|c|}{IIa}
								&	\multicolumn{2}{|c|}{Ia}		& 	\multicolumn{2}{|c|}{IIa}\\ 
\hline                  
\multicolumn{1}{|c|}{$\lambda$}	&	$|||e|||$			&	rate 	&	$|||e|||$			&	rate	
								&	$\|e\|_{L^2(\omg)}$	&	rate 	&	$\|e\|_{L^2(\omg)}$	&	rate	\\
\hline
\multicolumn{1}{|c|}{$2^3$}		&	1.11e-02 			& 	- 		&	3.37e-03				& 	-	
								&	1.12e-02 			& 	- 		&	3.47e-03				& 	-	\\
\hline
\multicolumn{1}{|c|}{$2^4$}		&	3.90e-03 			& 	1.51		&	1.19e-03				& 	1.50	
								&	3.92e-03 			& 	1.51		&	1.22e-03				& 	1.50	\\
\hline
\multicolumn{1}{|c|}{$2^5$}		&	1.37e-03 			&	1.51		&	4.19e-04  			& 	1.50		
								&	1.38e-03 			&	1.51		&	4.32e-04  			& 	1.50	\\
\hline
\multicolumn{1}{|c|}{$2^6$}		&   4.83e-04 			&	1.51		&	1.48e-04 			& 	1.51		
								&   4.87e-04 			&	1.51		&	1.52e-04 			& 	1.51		\\
\hline
\multicolumn{1}{|c|}{$2^7$}		&   1.69e-04 			&	1.52		&	5.16e-05 			& 	1.51		
								&   1.70e-04 			&	1.52		&	5.32e-05 			& 	1.51		\\
\hline
\end{tabular}
\caption{Dependence on the interaction radius $\lambda$ of the ``errors'' $|||e|||$ and $\|e\|_{L^2(\omgs)}$ and the rate of convergence of approximate nonlocal solutions. Results are provided for the test cases Ia and IIa.}
\label{tab:Arate}
\end{center}
\end{table}

\begin{table}[h!]
\begin{center}
\begin{tabular}{ c | c | c | c | c | c | c | c | c | }
\cline{2-9}
								&	\multicolumn{2}{|c|}{Ib}		& 	\multicolumn{2}{|c|}{IIb}
								&	\multicolumn{2}{|c|}{Ib}		& 	\multicolumn{2}{|c|}{IIb}\\ 
\hline                  
\multicolumn{1}{|c|}{$\lambda$}	&	$|||e|||$			&	rate 	&	$|||e|||$			&	rate	
								&	$\|e\|_{L^2(\omg)}$	&	rate 	&	$\|e\|_{L^2(\omg)}$	&	rate	\\
\hline
\multicolumn{1}{|c|}{$2^3$}		&	1.41e-01 			& 	- 		&	6.11e-02				& 	-	
								&	1.42e-01 			& 	- 		&	6.27e-02				& 	-	\\
\hline
\multicolumn{1}{|c|}{$2^4$}		&	7.58e-02 			& 	0.89		&	3.40e-02				& 	0.85		
								&	7.64e-02 			& 	0.90		&	3.49e-02				& 	0.84	\\
\hline
\multicolumn{1}{|c|}{$2^5$}		&	4.16e-02 			&	0.87		&	1.91e-02  			& 	0.83		
								&	4.19e-02 			&	0.87		&	1.95e-02  			& 	0.84	\\
\hline
\multicolumn{1}{|c|}{$2^6$}		&   2.29e-02 			&	0.86		&	1.07e-02 			& 	0.83		
								&   2.30e-02 			&	0.86		&	1.09e-02 			& 	0.84	\\
\hline
\multicolumn{1}{|c|}{$2^7$}		&   1.24e-02 			&	0.87		&	6.00e-03 			& 	0.83		
								&   1.25e-02 			&	0.88		&	6.06e-03 			& 	0.85	\\
\hline
\end{tabular}
\caption{Dependence on the interaction radius $\lambda$ of the ``errors'' $|||e|||$ and $\|e\|_{L^2(\omgs)}$ and the rate of convergence of approximate nonlocal solutions. Results are provided for the test cases Ib and IIb.}
\label{tab:Brate}
\end{center}
\end{table}

\section{Concluding remarks}
In this paper, we study a nonlocal diffusion operator which has as a special case the fractional Laplacian operator $(-\Delta)^s$; exploiting a nonlocal vector calculus, we show that the solutions of nonlocal problems for the operator $\mcL$ converge to the solutions of fractional Laplacian problems as the nonlocal interactions become infinite. Through several numerical examples, we compare the nonlocal solutions with the solution of the fractional Laplacian equation on bounded domains, illustrate the theoretical results, and show that solving nonlocal problems is a viable alternative to solving fractional differential problems which feature infinite-volume constraints. 

A possible extension of this work consists in improving the rate of convergence, with respect to the grid size, of finite element approximations. The analytic solution of the fractional Laplacian equation (see e.g. \eqref{eq:intu1}) presents an abrupt change at the end points of the domain $\omgs$; in fact, it has a singular derivative in $x=-1,1$. Incorporating the form of the singularity in the numerical scheme would yield greatly enhanced convergence; among the possible approaches we mention the {\it singular basis function method} \cite{geor:89,olge:91}. In this approach, a set of supplementary functions that reproduce the functional form, assumed known, of the solution singularity is added to the ordinary finite element basis functions. Given the asymptotic expansion of the analytic solution as $x\to -1,1$, $u(x)=\sum_{k=1}^\infty a_k x^{\lambda_k}$, where $a_k$ and $\lambda_k$ are the singular coefficients and exponents respectively, the supplementary basis functions have the general form 
$$\psi_k(x)=b(x)x^{\lambda_k},$$ 
$b(x)$ being an optional blending function.

Another possible follow-up of the present work is to treat the nonlocal time-dependent diffusion equation  
\begin{equation}\label{eq:unsteady}
\left\{\begin{array}{ll}
u_t = \mcL u &  {\rm in} \; \omgs \\
[3mm]
u(\cdot,t) = 0 & {\rm in} \; \mbR^n\backslash\omgs\\
[3mm]
u(\cdot,0) = u_0 &  {\rm in} \; \omgs.
\end{array}\right.
\end{equation}
This is a problem of interest in applications involving jump processes \cite{burc:11}. In fact, \eqref{eq:unsteady} is the {\it master equation} for a jump process, i.e., the deterministic equation that determines the time evolution of the probability density function restricted to a bounded domain $\omgs$. As done for the steady case, the solution $u$ of \eqref{eq:unsteady} can be approximated by the solution $\wtu$ of 
\begin{equation}\label{eq:unsteady_nl}
\left\{\begin{array}{ll}
\wtu_t = \widetilde\mcL \wtu &  {\rm in} \; \omgs \\
[3mm]
\wtu(\cdot,t) = 0 & {\rm in} \; \omgc\\
[3mm]
\wtu(\cdot,0) = u_0 &  {\rm in} \; \omgs,
\end{array}\right.
\end{equation}
where $\widetilde\mcL$ is defined as in \S\ref{nonlocal_fractional}. A fractional time derivative might be considered as well.

Although in this work fractional Laplacian problems posed on all of $\mbR^n$ are not considered, such problems are indeed of interest so that the role of nonlocal diffusion operators for such problems would be interesting to explore. In this case, one would want to study the behavior of solutions of such problems as the domain $\Omega$, which for computational purposes has to be chosen to be finite, increases in size.

We have only considered continuous Galerkin discretizations of the nonlocal diffusion problem, e.g., we have used continuous piecewise-linear finite element bases. For $s\le 1/2$, {\em discontinuous} Galerkin methods are also conforming, i.e., satisfy \eqref{conf}. Thus it would be of interest to study such discretizations, especially in view of the fact that for $s\le 1/2$, both  nonlocal diffusion and fractional Laplacian problems admit solutions containing jump discontinuities.

Even if preliminary in nature, our numerical results presented here suffice to illustrate the theory and show the viability of using volume-constrained nonlocal diffusion problems as a means for defining approximations of fractional Laplacian problems posed on finite domains. Obviously, a natural follow-up of this work is to extend the numerical simulations to two-dimensional and three-dimensional settings for which computations would be even more challenging.

\end{document}